\documentclass[11pt,reqno,letterpaper]{amsart}
\usepackage{color}
\usepackage[colorlinks=true, allcolors=blue,backref=page]{hyperref}
\usepackage{amsmath, amssymb, amsthm}
\usepackage{mathrsfs}
\usepackage{mathtools}
\usepackage[noabbrev,capitalize,nameinlink]{cleveref}
\crefname{equation}{}{}
\usepackage{fullpage}
\usepackage[noadjust]{cite}
\usepackage{graphics}
\usepackage{pifont}

\usepackage{tikz}
\usetikzlibrary{calc}
\usetikzlibrary{shapes.geometric}

\usepackage{bbm}
\usepackage{multicol} % Include multicol package for column alignment
\usepackage[T1]{fontenc}

\usepackage{environ}
\usepackage{framed}
\usepackage{url}
\usepackage[linesnumbered,ruled,vlined]{algorithm2e}
\usepackage[noend]{algpseudocode}
\usepackage[labelfont=bf]{caption}
\usepackage[framemethod=tikz]{mdframed}
\usepackage{appendix}
\usepackage{graphicx}
\usepackage[textsize=tiny]{todonotes}
\usepackage{tcolorbox}
\usepackage{enumerate}
\allowdisplaybreaks[1]
\usepackage{stmaryrd}
\usepackage{colonequals}
\usepackage[margin=1in]{geometry}
\usepackage{caption}
\usepackage{subcaption}

\usepackage[shortlabels]{enumitem}
\crefformat{enumi}{#2#1#3}
\crefrangeformat{enumi}{#3#1#4 to~#5#2#6}
\crefmultiformat{enumi}{#2#1#3}
{ and~#2#1#3}{, #2#1#3}{ and~#2#1#3}

\DeclareSymbolFont{symbolsC}{U}{pxsyc}{m}{n}
\SetSymbolFont{symbolsC}{bold}{U}{pxsyc}{bx}{n}
\DeclareFontSubstitution{U}{pxsyc}{m}{n}
\DeclareMathSymbol{\medcircle}{\mathbin}{symbolsC}{7}

\crefname{algocf}{Algorithm}{Algorithms}

\crefname{equation}{}{} %remove ``Equation''
\AtBeginEnvironment{appendices}{\crefalias{section}{appendix}} %appendices

\usepackage[color,final]{showkeys} %add in 'final' into parameter to remove showkeys

\colorlet{refkey}{orange!20}
\colorlet{labelkey}{blue!30}

\crefname{algocf}{Algorithm}{Algorithms}

\numberwithin{equation}{section}
\newtheorem{theorem}{Theorem}[section]

\newtheorem{lemma}[theorem]{Lemma}
\newtheorem{claim}[theorem]{Claim}

\crefname{claim}{Claim}{Claims}

\newtheorem{corollary}[theorem]{Corollary}
\newtheorem{conjecture}[theorem]{Conjecture}
\newtheorem*{question*}{Question}

\theoremstyle{definition}
\newtheorem{definition}[theorem]{Definition}

\newtheorem*{definition*}{Definition}

\theoremstyle{remark}
\newtheorem{remark}[theorem]{Remark}

\setlist[enumerate,1]{label={(\arabic*)}}

\newlist{enumthm}{enumerate}{1}
\setlist[enumthm]{label=\textup{(\roman*)},ref=\thethm(\roman*)}
\Crefname{enumthmi}{Theorem}{Theorems}

\newcommand{\norm}[1]{\left\lVert#1\right\rVert}

\newcommand{\R}{\mathbb R}

\newcommand{\N}{\mathbb N}

\newcommand{\mb}{\mathbb}

\newcommand{\mc}{\mathcal}

\newcommand{\one}{{\bf 1}}

\newcommand{\eps}{\varepsilon}
\newcommand{\codeg}{\mathrm{codeg}}

\let\originalleft\left
\let\originalright\right
\renewcommand{\left}{\mathopen{}\mathclose\bgroup\originalleft}
\renewcommand{\right}{\aftergroup\egroup\originalright}

\allowdisplaybreaks

\newcounter{GProperties}

\newcommand{\bowtieFigure}{
\begin{tikzpicture}
\tikzset{every node/.style={fill=black, circle, inner sep=2pt}}

\pgfmathsetmacro{\L}{3} % side length
\pgfmathsetmacro{\dx}{\L/2} % horizontal offset from u and v
\pgfmathsetmacro{\h}{sqrt(\L)/2} % height of equilateral triangle

\coordinate (u) at (-\L, 0);
\coordinate (a1) at (-\dx, \h);
\coordinate (a2) at (-\dx, -\h);

\coordinate (v) at (\L, 0);
\coordinate (b1) at (\dx, \h);
\coordinate (b2) at (\dx, -\h);

\coordinate (c) at (0, 0);

\begin{scope}[fill=gray!50]
    \fill (u) -- (a2) -- (a1) -- cycle;
    \fill (c) -- (b1) -- (b2) -- cycle;
\end{scope}

\draw (u) -- (a1);
\draw (u) -- (a2);
\draw (a1) -- (a2);
\draw (a1) -- (c);
\draw (a2) -- (c);
\draw (v) -- (b1);
\draw (v) -- (b2);
\draw (b1) -- (b2);
\draw (b1) -- (c);
\draw (b2) -- (c);

\node at (u) [label=$u$] {};
\node at (v) [label=$v$] {};
\node at (a1) [label=$a_1$] {};
\node at (a2) [label=below:$a_2$] {};
\node at (b1) [label=$b_1$] {};
\node at (b2) [label=below:$b_2$] {};
\node at (c) [label=$c$] {};

\end{tikzpicture}
}

\newcommand{\octagon}{
\begin{tikzpicture}[scale=2]
  \tikzset{every node/.style={fill=black, circle, inner sep=2pt}}

  \foreach \i/\name in {247.5/v, 292.5/u, 337.5/w1, 22.5/w2, 67.5/w3, 112.5/w4, 157.5/w5, 202.5/w6} {
    \coordinate (\name) at ({cos(\i)}, {sin(\i)});
  }
  \coordinate (c) at (0,0);

  \foreach \i/\j in {u/v,w1/w2,w3/w4,w5/w6} {
    \fill[fill=gray!50] (\i) -- (\j) -- (c) -- cycle;
  }

  \foreach \name in {v,u,w1,w2,w3,w4,w5,w6} {
    \draw (c) -- (\name);
  }

  \draw (v) -- (u) -- (w1) -- (w2) -- (w3) -- (w4) -- (w5) -- (w6) -- (v);

  \foreach \name/\nameLabel in {c/$c$,v/$v$,u/$u$,w1/$w_1$,w2/$w_2$,w3/$w_3$,w4/$w_4$,w5/$w_5$,w6/$w_6$} {
    \node at (\name) [label=\nameLabel] {};
  }
  
\end{tikzpicture}}

\newcommand{\triangleAndIso}{
\begin{tikzpicture}[baseline]
\fill[black] (0,0) circle [radius=1pt];
\fill[black] (.2,0) circle [radius=1pt];
\fill[black] (0,.2) circle [radius=1pt];
\fill[black] (.2,.2) circle [radius=1pt];
\draw (0,0) -- (.2,0) -- (0,.2) -- (0,0);
\end{tikzpicture}
}

\newcommand{\triangleWithTail}{
\begin{tikzpicture}[baseline]
\fill[black] (0,0) circle [radius=1pt];
\fill[black] (.2,0) circle [radius=1pt];
\fill[black] (0,.2) circle [radius=1pt];
\fill[black] (.2,.2) circle [radius=1pt];
\draw (0,.2) -- (.2,0) -- (0,0) -- (0,.2) -- (.2,.2);
\end{tikzpicture}
}

\newcommand{\diamondGraph}{
\begin{tikzpicture}[baseline]
\fill[black] (0,0) circle [radius=1pt];
\fill[black] (.2,0) circle [radius=1pt];
\fill[black] (0,.2) circle [radius=1pt];
\fill[black] (.2,.2) circle [radius=1pt];
\draw (0,0) -- (.2,0) -- (.2,.2) -- (0,.2) -- (0,0) -- (.2,.2);
\end{tikzpicture}
}

\newcommand{\wThree}{
\begin{tikzpicture}[scale=2.5]
  \tikzset{every node/.style={fill=black, circle, inner sep=2pt}}

  \foreach \i/\name in {247.5/w7, 292.5/w0, 337.5/w1, 22.5/w2, 67.5/w3, 112.5/w4, 157.5/w5, 202.5/w6} {
    \coordinate (\name) at (-{cos(\i)}, {sin(\i)});
  }
  \coordinate (c) at (0,0);
  \coordinate (d) at ({2*cos(22.5)},-{cos(22.5)});

  \foreach \sign in {-1,1} {
    \pgfmathsetmacro{\xOffset}{\sign * cos(22.5)}
    \coordinate (offset) at (\xOffset, 0);

    \foreach \name in {w7,w0,w1,w2,w3,w4,w5,w6} {
      \draw ($(offset) + (c)$) -- ($(offset) + (\name)$);
    }

    \draw 
      ($(offset)+(w7)$) -- ($(offset)+(w0)$) -- ($(offset)+(w1)$) -- 
      ($(offset)+(w2)$) -- ($(offset)+(w3)$) -- ($(offset)+(w4)$) -- 
      ($(offset)+(w5)$) -- ($(offset)+(w6)$) -- cycle;
  }

  \coordinate (offset) at ({cos(22.5)},0);

  \draw ($(offset) + (w6)$) -- (d) -- ($(offset) + (w7)$);
  
  \foreach \name/\labeltext in {c/$c_a$,w3/$a_4$,w4/$a_5$,w5/$a_6$,w6/$a_7$,w7/$a_0$,w0/$a_1$} {
    \node at ($(offset)+(\name)$) [label=\labeltext] {};
  }
  \node at ($(offset)+(w1)$) [label=below:{$a_2=b_0$}] {};
  \node at ($(offset)+(w2)$) [label={$a_3=b_7$}] {};
  
  \coordinate (offset) at (-{cos(22.5)},0);
  \foreach \name/\labeltext in {c/$c_b$,w3/$b_5$,w4/$b_6$,w7/$b_1$,w0/$b_2$,w1/$b_3$,w2/$b_4$} {
    \node at ($(offset)+(\name)$) [label=\labeltext] {};
  }

  \node at (d) [label=right:$d$] {};
\end{tikzpicture}
}

\newcommand{\WTwoSeven}{
\begin{tikzpicture}[scale=1.2]
    \node (v) at (2,3) [circle, draw, fill=red!50, minimum size=10pt] {$d$};
    \node (u) at (0,-1) [regular polygon, regular polygon sides=4, draw,fill=blue!50, minimum size=10pt] {$a_5$};
    \node (w_i) at (0,1) [regular polygon, regular polygon sides=4, draw,fill=blue!50, minimum size=10pt] {$a_6$};
    \node (c_1) at (2,0) [regular polygon, regular polygon sides=5, draw,fill=gray!50, minimum size=10pt] {$c_a$};
    \node (a_1) at (1,2) [circle, draw,fill=red!50, minimum size=5pt] {$a_7$};
    \node (a_2) at (3,2) [regular polygon, regular polygon sides=5, draw,fill=gray!50, minimum size=10pt] {$a_0$};
    \node (a_3) at (4,1) [regular polygon, regular polygon sides=5, draw,fill=gray!50, minimum size=10pt] {$a_1/b_0$};
    \node (a_4) at (4,-1) [circle, draw,fill=red!50, minimum size=5pt] {$a_2/b_7$};
    \node (a_5) at (3,-2) [regular polygon, regular polygon sides=4, draw,fill=blue!50, minimum size=10pt] {$a_3$};
    \node (a_6) at (1,-2) [regular polygon, regular polygon sides=4, draw,fill=blue!50, minimum size=10pt] {$a_4$};
    \node (c_2) at (6,0) [regular polygon, regular polygon sides=5, draw,fill=gray!50, minimum size=10pt] {$c_b$};
    \node (b_1) at (5,2) [regular polygon, regular polygon sides=4, draw,fill=blue!50, minimum size=10pt] {$b_1$};
    \node (b_2) at (7,2) [regular polygon, regular polygon sides=4, draw,fill=blue!50, minimum size=10pt] {$b_2$};
    \node (b_3) at (8,1) [regular polygon, regular polygon sides=4, draw,fill=blue!50, minimum size=10pt] {$b_3$};
    \node (b_4) at (8,-1) [regular polygon, regular polygon sides=4, draw,fill=blue!50, minimum size=10pt] {$b_4$};
    \node (b_5) at (7,-2) [regular polygon, regular polygon sides=4, draw,fill=blue!50, minimum size=10pt] {$b_5$};
    \node (b_6) at (5,-2) [circle, draw,fill=red!50, minimum size=5pt] {$b_6$};
    \draw[thick, blue] (w_i) -- (u);
    \draw (v) -- (a_2);
    \draw[thick, red] (a_1) -- (v);
    \draw[thick, blue] (c_1) -- (u);
    \draw[thick, blue](c_1) -- (w_i);
    \draw (c_1) -- (a_1);
    \draw (c_1) -- (a_2);
    \draw (c_1) -- (a_3);
    \draw (c_1) -- (a_4);
    \draw[thick, blue] (c_1) -- (a_5);
    \draw[thick, blue] (c_1) -- (a_6);
    \draw[thick, blue] (w_i) -- (a_1);
    \draw (a_1) -- (a_2);
    \draw (a_2) -- (a_3);
    \draw (a_3) -- (a_4);
    \draw[thick, blue] (a_4) -- (a_5);
    \draw[thick, blue] (a_5) -- (a_6);
    \draw[thick, blue] (a_6) -- (u);
    \draw (c_2) -- (a_4);
    \draw (c_2) -- (a_3);
    \draw[thick, blue] (c_2) -- (b_1);
    \draw[thick, blue] (c_2) -- (b_2);
    \draw[thick, blue] (c_2) -- (b_3);
    \draw[thick, blue] (c_2) -- (b_4);
    \draw[thick, blue] (c_2) -- (b_5);
    \draw (c_2) -- (b_6);
    \draw[thick, blue] (a_3) -- (b_1);
    \draw[thick, blue] (b_1) -- (b_2);
    \draw[thick, blue] (b_2) -- (b_3);
    \draw[thick, blue] (b_3) -- (b_4);
    \draw[thick, blue] (b_4) -- (b_5);
    \draw[thick, blue] (b_5) -- (b_6);
    \draw[thick, red](b_6) -- (a_4);   
\end{tikzpicture}
}

\newcommand{\deltaIIllustration}{
\begin{tikzpicture}[scale=1.2]
    \node (u) at (1,0) [circle, draw, fill=red!50, minimum size=10pt] {$v/w_7$};
    \node (w_i) at (-1,0) [circle, draw,fill=red!50, minimum size=10pt] {$u/w_0$};
    \node (c_1) at (0,2) [regular polygon, regular polygon sides=4, shape border rotate=45, draw,fill=green!50, minimum size=5pt] {$c$};
    \node (a_1) at (-2,1)[regular polygon, regular polygon sides=4, shape border rotate=45, draw,fill=green!50, minimum size=5pt] {$w_1$};
    \node (a_2) at (-2,3) [regular polygon, regular polygon sides=4, shape border rotate=45, draw,fill=green!50, minimum size=5pt] {$w_2$};
    \node (a_3) at (-1,4) [star, star points=5, draw, fill=orange!50, minimum size=5pt] {$w_3$};
    \node (a_4) at (1,4) [regular polygon, regular polygon sides=4, draw,fill=blue!50, minimum size=10pt] {$w_4$};
    \node (a_5) at (2,3) [regular polygon, regular polygon sides=4, draw,fill=blue!50, minimum size=10pt] {$w_5$};
    \node (a_6) at (2,1) [regular polygon, regular polygon sides=4, draw,fill=blue!50, minimum size=10pt] {$w_6$};   
    \draw[thick, red] (w_i) -- (u);
    \draw[thick, green] (c_1) -- (u);
    \draw[thick, green] (c_1) -- (w_i);
    \draw[thick, green] (c_1) -- (a_1);
    \draw[thick, green] (c_1) -- (a_2);
    \draw[thick, orange]  (c_1) -- (a_3);
    \draw[thick, blue]  (c_1) -- (a_4);
    \draw[thick, blue]  (c_1) -- (a_5);
    \draw[thick, blue]  (c_1) -- (a_6);
    \draw[thick, green] (w_i) -- (a_1);
    \draw[thick, green] (a_1) -- (a_2);
    \draw[thick, orange] (a_2) -- (a_3);
    \draw[thick, blue] (a_3) -- (a_4);
    \draw[thick, blue]  (a_4) -- (a_5);
    \draw[thick, blue]  (a_5) -- (a_6);
    \draw[thick, blue] (a_6) -- (u);
\end{tikzpicture}
}

\newcommand{\WOneTwo}{
\begin{tikzpicture}[scale=1.4]
    \node (v) at (-4,3) [circle, draw, fill=red!50, minimum size=10pt] {$d$};
    \node (u) at (1,0) [regular polygon, regular polygon sides=4, draw,fill=blue!50, minimum size=10pt] {$a_4$};
    \node (w_i) at (-1,0) [regular polygon, regular polygon sides=4, draw,fill=blue!50, minimum size=10pt] {$a_5$};
    \node (c_1) at (0,2) [regular polygon, regular polygon sides=5, draw,fill=gray!50, minimum size=10pt] {$c_a$};
    \node (a_1) at (-2,1) [regular polygon, regular polygon sides=4, draw,fill=blue!50, minimum size=10pt] {$a_6$};
    \node (a_2) at (-2,3) [circle, draw,fill=red!50, minimum size=5pt] {$a_7$};
    \node (a_3) at (-1,4) [regular polygon, regular polygon sides=5, draw,fill=gray!50, minimum size=10pt] {$a_0/b_0$};
    \node (a_4) at (1,4) [regular polygon, regular polygon sides=5, draw,fill=gray!50, minimum size=10pt] {$a_1/b_7$};
    \node (a_5) at (2,3) [regular polygon, regular polygon sides=4, draw,fill=blue!50, minimum size=10pt] {$a_2$};
    \node (a_6) at (2,1) [regular polygon, regular polygon sides=4, draw,fill=blue!50, minimum size=10pt] {$a_3$};
    \node (c_2) at (0,6) [regular polygon, regular polygon sides=5, draw,fill=gray!50, minimum size=10pt] {$c_b$};
    \node (b_1) at (-2,5) [regular polygon, regular polygon sides=5, draw,fill=gray!50, minimum size=10pt] {$b_1/\vartheta_{0}$};
    \node (b_2) at (-2,7) [regular polygon, regular polygon sides=5, draw,fill=gray!50, minimum size=10pt] {$b_2/\vartheta_{7}$};
    \node (b_3) at (-1,8) [regular polygon, regular polygon sides=4, draw,fill=blue!50, minimum size=10pt] {$b_3$};
    \node (b_4) at (1,8) [regular polygon, regular polygon sides=4, draw,fill=blue!50, minimum size=10pt] {$b_4$};
    \node (b_5) at (2,7) [regular polygon, regular polygon sides=4, draw,fill=blue!50, minimum size=10pt] {$b_5$};
    \node (b_6) at (2,5) [regular polygon, regular polygon sides=4, draw,fill=blue!50, minimum size=10pt] {$b_6$};
    \node (c_3) at (-4,6) [regular polygon, regular polygon sides=5, draw,fill=gray!50, minimum size=10pt] {$c_{\vartheta}$};
    \node (d_1) at (-3,4) [regular polygon, regular polygon sides=4, draw,fill=blue!50, minimum size=10pt] {$\vartheta_{1}$};
    \node (d_2) at (-5,4) [regular polygon, regular polygon sides=4, draw,fill=blue!50, minimum size=10pt] {$\vartheta_{2}$};
    \node (d_3) at (-6,5) [regular polygon, regular polygon sides=4, draw,fill=blue!50, minimum size=10pt] {$\vartheta_{3}$};
    \node (d_4) at (-6,7) [regular polygon, regular polygon sides=4, draw,fill=blue!50, minimum size=10pt] {$\vartheta_{4}$};
    \node (d_5) at (-5,8) [circle, draw,fill=red!50, minimum size=5pt] {$\vartheta_{5}$};
    \node (d_6) at (-3,8) [circle, draw,fill=red!50, minimum size=5pt] {$\vartheta_{6}$};
    \draw[thick, blue] (w_i) -- (u);
    \draw (v) -- (a_3);
    \draw[thick, red] (a_2) -- (v);
    \draw[thick, blue] (c_1) -- (u);
    \draw[thick, blue](c_1) -- (w_i);
    \draw[thick, blue] (c_1) -- (a_1);
    \draw (c_1) -- (a_2);
    \draw (c_1) -- (a_3);
    \draw (c_1) -- (a_4);
    \draw[thick, blue] (c_1) -- (a_5);
    \draw[thick, blue] (c_1) -- (a_6);
    \draw[thick, blue] (w_i) -- (a_1);
    \draw[thick, blue] (a_1) -- (a_2);
    \draw (a_2) -- (a_3);
    \draw (a_3) -- (a_4);
    \draw[thick, blue] (a_4) -- (a_5);
    \draw[thick, blue] (a_5) -- (a_6);
    \draw[thick, blue] (a_6) -- (u);
    \draw (c_2) -- (a_4);
    \draw (c_2) -- (a_3);
    \draw (c_2) -- (b_1);
    \draw (c_2) -- (b_2);
    \draw[thick, blue] (c_2) -- (b_3);
    \draw[thick, blue] (c_2) -- (b_4);
    \draw[thick, blue] (c_2) -- (b_5);
    \draw[thick, blue] (c_2) -- (b_6);
    \draw (a_3) -- (b_1);
    \draw (b_1) -- (b_2);
    \draw[thick, blue] (b_2) -- (b_3);
    \draw[thick, blue] (b_3) -- (b_4);
    \draw[thick, blue] (b_4) -- (b_5);
    \draw[thick, blue] (b_5) -- (b_6);
    \draw[thick, blue] (b_6) -- (a_4);
    \draw (c_3) -- (b_1);
    \draw (c_3) -- (b_2);
    \draw[thick, blue] (c_3) -- (d_1);
    \draw[thick, blue] (c_3) -- (d_2);
    \draw[thick, blue] (c_3) -- (d_3);
    \draw[thick, blue] (c_3) -- (d_4);
    \draw (c_3) -- (d_5);
    \draw (c_3) -- (d_6);
    \draw[thick, blue] (b_1) -- (d_1);
    \draw[thick, blue] (d_1) -- (d_2);
    \draw[thick, blue] (d_2) -- (d_3);
    \draw[thick, blue] (d_3) -- (d_4);
    \draw[thick, blue] (d_4) -- (d_5);
    \draw[thick, red] (d_5) -- (d_6);
    \draw (d_6) -- (b_2);
    
\end{tikzpicture}
}

\newcommand{\HOne}{
\begin{tikzpicture}[scale=1.4]
    \node (v) at (-3,2) [circle, draw, fill=red!50, minimum size=10pt] {$v$};
    \node (u) at (1,0) [regular polygon, regular polygon sides=4, draw,fill=blue!50, minimum size=10pt] {$a_5$};
    \node (w_i) at (-1,0) [regular polygon, regular polygon sides=4, draw,fill=blue!50, minimum size=10pt] {$a_6$};
    \node (c_1) at (0,2) [regular polygon, regular polygon sides=4, shape border rotate=45, draw,fill=green!50, minimum size=5pt] {$2,c_a$};
    \node (a_1) at (-2,1) [circle, draw,fill=red!50, minimum size=5pt] {$u,a_7$};
    \node (a_2) at (-2,2.5) [regular polygon, regular polygon sides=4, shape border rotate=45, draw,fill=green!50, minimum size=5pt] {$1,a_0$};
    \node (a_3) at (-1,4) [regular polygon, regular polygon sides=4, shape border rotate=45, draw,fill=green!50, minimum size=5pt] {$3,\varsigma,a_1,b_0$};
    \node (a_4) at (1.2,4) [regular polygon, regular polygon sides=6, draw,fill=purple!50, minimum size=5pt] {$a_2,b_7$};
    \node (a_5) at (2,3) [regular polygon, regular polygon sides=4, draw,fill=blue!50, minimum size=5pt] {$a_3$};
    \node (a_6) at (2,1) [regular polygon, regular polygon sides=4, draw,fill=blue!50, minimum size=5pt] {$a_4$};
    \node (c_2) at (0,6) [regular polygon, regular polygon sides=6, draw,fill=purple!50, minimum size=5pt] {$c_b$};
    \node (b_1) at (-2,5.7) [star, star points=5, draw, fill=orange!50, minimum size=5pt] {$\beta,b_1 $};
    \node (b_2) at (-2,7) [regular polygon, regular polygon sides=4, draw,fill=blue!50, minimum size=5pt] {$b_2$};
    \node (b_3) at (-1,8) [regular polygon, regular polygon sides=4, draw,fill=blue!50, minimum size=5pt] {$b_3$};
    \node (b_4) at (1,8) [regular polygon, regular polygon sides=4, draw,fill=blue!50, minimum size=5pt] {$b_4$};
    \node (b_5) at (2,7) [regular polygon, regular polygon sides=4, draw,fill=blue!50, minimum size=5pt] {$b_5$};
    \node (b_6) at (2,5) [regular polygon, regular polygon sides=4, draw,fill=blue!50, minimum size=5pt] {$b_6$};    
    \draw[thick, blue] (w_i) -- (u);
    \draw[thick, green] (v) -- (a_2);
    \draw[thick, red] (a_1) -- (v);
    \draw[thick, blue] (c_1) -- (u);
    \draw[thick, blue](c_1) -- (w_i);
    \draw[thick, green] (c_1) -- (a_1);
    \draw[thick, green] (c_1) -- (a_2);
    \draw[thick, green] (c_1) -- (a_3);
    \draw[thick, purple] (c_1) -- (a_4);
    \draw[thick, blue] (c_1) -- (a_5);
    \draw[thick, blue] (c_1) -- (a_6);
    \draw[thick, blue] (w_i) -- (a_1);
    \draw[thick, green] (a_1) -- (a_2);
    \draw[thick, green] (a_2) -- (a_3);
    \draw[thick, purple] (a_3) -- (a_4);
    \draw[thick, blue] (a_4) -- (a_5);
    \draw[thick, blue] (a_5) -- (a_6);
    \draw[thick, blue] (a_6) -- (u);
    \draw[thick, purple] (c_2) -- (a_4);
    \draw[thick, purple] (c_2) -- (a_3);
    \draw[thick, purple] (c_2) -- (b_1);
    \draw[thick, blue] (c_2) -- (b_2);
    \draw[thick, blue] (c_2) -- (b_3);
    \draw[thick, blue] (c_2) -- (b_4);
    \draw[thick, blue] (c_2) -- (b_5);
    \draw[thick, blue] (c_2) -- (b_6);
    \draw[thick, orange] (a_3) -- (b_1);
    \draw[thick, blue] (b_1) -- (b_2);
    \draw[thick, blue] (b_2) -- (b_3);
    \draw[thick, blue] (b_3) -- (b_4);
    \draw[thick, blue] (b_4) -- (b_5);
    \draw[thick, blue] (b_5) -- (b_6);
    \draw[thick, blue] (b_6) -- (a_4);   
\end{tikzpicture}
}

\newcommand{\HTwo}{
\begin{tikzpicture}[scale=1.4]
    \node (v) at (3,2) [circle, draw, fill=red!50, minimum size=10pt] {$v$};
    \node (u) at (1,0) [regular polygon, regular polygon sides=4, draw,fill=blue!50, minimum size=10pt] {$a_1$};
    \node (w_i) at (-1,0) [regular polygon, regular polygon sides=4, draw,fill=blue!50, minimum size=10pt] {$a_2$};
    \node (c_1) at (0,2) [regular polygon, regular polygon sides=4, shape border rotate=45, draw,fill=green!50, minimum size=5pt] {$2,c_a$};
    \node (a_1) at (-2,1) [regular polygon, regular polygon sides=4, draw,fill=blue!50, minimum size=5pt] {$a_3$};
    \node (a_2) at (-2,3) [regular polygon, regular polygon sides=4, draw,fill=blue!50, minimum size=5pt] {$a_4$};
    \node (a_3) at (-1.2,4) [regular polygon, regular polygon sides=6, draw,fill=purple!50, minimum size=5pt] {$a_5,b_0$};
    \node (a_4) at (1,4) [regular polygon, regular polygon sides=4, shape border rotate=45, draw,fill=green!50, minimum size=5pt] {$3,\varsigma,a_6,b_7$};
    \node (a_5) at (2,2.7) [circle, draw,fill=red!50, minimum size=5pt] {$u,a_7$};
    \node (a_6) at (2,1) [regular polygon, regular polygon sides=4, shape border rotate=45, draw,fill=green!50, minimum size=5pt] {$1,a_0$};
    \node (c_2) at (0,6) [regular polygon, regular polygon sides=6, draw,fill=purple!50, minimum size=5pt] {$c_b$};
    \node (b_1) at (-2,5) [regular polygon, regular polygon sides=4, draw,fill=blue!50, minimum size=5pt] {$b_1$};
    \node (b_2) at (-2,7) [regular polygon, regular polygon sides=4, draw,fill=blue!50, minimum size=5pt] {$b_2$};
    \node (b_3) at (-1,8) [regular polygon, regular polygon sides=4, draw,fill=blue!50, minimum size=5pt] {$b_3$};
    \node (b_4) at (1,8) [regular polygon, regular polygon sides=4, draw,fill=blue!50, minimum size=5pt] {$b_4$};
    \node (b_5) at (2,7) [regular polygon, regular polygon sides=4, draw,fill=blue!50, minimum size=5pt] {$b_5$};
    \node (b_6) at (2,5.7) [star, star points=5, draw, fill=orange!50, minimum size=5pt] {$\beta,b_6$};    
    \draw[thick, blue] (w_i) -- (u);
    \draw[thick, red] (v) -- (a_5);
    \draw[thick, green] (a_6) -- (v);
    \draw[thick, blue] (c_1) -- (u);
    \draw[thick, blue] (c_1) -- (w_i);
    \draw[thick, blue] (c_1) -- (a_1);
    \draw[thick, blue] (c_1) -- (a_2);
    \draw[thick, purple]  (c_1) -- (a_3);
    \draw[thick, green]  (c_1) -- (a_4);
    \draw[thick, green]  (c_1) -- (a_5);
    \draw[thick, green]  (c_1) -- (a_6);
    \draw[thick, blue] (w_i) -- (a_1);
    \draw[thick, blue] (a_1) -- (a_2);
    \draw[thick, blue] (a_2) -- (a_3);
    \draw[thick, purple] (a_3) -- (a_4);
    \draw[thick, green]  (a_4) -- (a_5);
    \draw[thick, green]  (a_5) -- (a_6);
    \draw[thick, blue] (a_6) -- (u);
    \draw[thick, purple] (c_2) -- (a_4);
    \draw[thick, purple] (c_2) -- (a_3);
    \draw[thick, blue] (c_2) -- (b_1);
    \draw[thick, blue] (c_2) -- (b_2);
    \draw[thick, blue] (c_2) -- (b_3);
    \draw[thick, blue] (c_2) -- (b_4);
    \draw[thick, blue] (c_2) -- (b_5);
    \draw[thick, purple] (c_2) -- (b_6);
    \draw[thick, blue] (a_3) -- (b_1);
    \draw[thick, blue] (b_1) -- (b_2);
    \draw[thick, blue] (b_2) -- (b_3);
    \draw[thick, blue] (b_3) -- (b_4);
    \draw[thick, blue] (b_4) -- (b_5);
    \draw[thick, blue] (b_5) -- (b_6);
    \draw[thick, orange](b_6) -- (a_4);   
\end{tikzpicture}
}

\newcommand{\HThree}{
\begin{tikzpicture}[scale=1.4]
    \node (v) at (-3,4) [circle, draw, fill=red!50, minimum size=10pt] {$v$};
    \node (u) at (1,0) [regular polygon, regular polygon sides=4, draw,fill=blue!50, minimum size=10pt] {$a_4$};
    \node (w_i) at (-1,0) [regular polygon, regular polygon sides=4, draw,fill=blue!50, minimum size=10pt] {$a_5$};
    \node (c_1) at (0,2) [regular polygon, regular polygon sides=4, shape border rotate=45, draw,fill=green!50, minimum size=5pt] {$2,c_a$};
    \node (a_1) at (-2,1) [regular polygon, regular polygon sides=4, draw,fill=blue!50, minimum size=5pt] {$a_6$};
    \node (a_2) at (-2,2.8) [circle, draw,fill=red!50, minimum size=5pt] {$u,a_7$};
    \node (a_3) at (-1,4) [regular polygon, regular polygon sides=4, shape border rotate=45, draw,fill=green!50, minimum size=5pt] {$1,\varsigma,a_0,b_0$};
    \node (a_4) at (1.2,4) [regular polygon, regular polygon sides=6, draw,fill=purple!50, minimum size=5pt] {$a_1,b_7$};
    \node (a_5) at (2,3) [regular polygon, regular polygon sides=4, draw,fill=blue!50, minimum size=5pt] {$a_2$};
    \node (a_6) at (2,1) [regular polygon, regular polygon sides=4, draw,fill=blue!50, minimum size=5pt] {$a_3$};
    \node (c_2) at (0,6) [regular polygon, regular polygon sides=6, draw,fill=purple!50, minimum size=5pt] {$c_b$};
    \node (b_1) at (-2,5.7) [star, star points=5, draw, fill=orange!50, minimum size=5pt] {$\beta,b_1$};
    \node (b_2) at (-2,7) [regular polygon, regular polygon sides=4, draw,fill=blue!50, minimum size=5pt] {$b_2$};
    \node (b_3) at (-1,8) [regular polygon, regular polygon sides=4, draw,fill=blue!50, minimum size=5pt] {$b_3$};
    \node (b_4) at (1,8) [regular polygon, regular polygon sides=4, draw,fill=blue!50, minimum size=5pt] {$b_4$};
    \node (b_5) at (2,7) [regular polygon, regular polygon sides=4, draw,fill=blue!50, minimum size=5pt] {$b_5$};
    \node (b_6) at (2,5) [regular polygon, regular polygon sides=4, draw,fill=blue!50, minimum size=5pt] {$b_6$};    
    \draw[thick, blue] (w_i) -- (u);
    \draw[thick, green] (v) -- (a_3);
    \draw[thick, red] (a_2) -- (v);
    \draw[thick, blue] (c_1) -- (u);
    \draw[thick, blue](c_1) -- (w_i);
    \draw[thick, blue] (c_1) -- (a_1);
    \draw[thick, green] (c_1) -- (a_2);
    \draw[thick, green] (c_1) -- (a_3);
    \draw[thick, purple] (c_1) -- (a_4);
    \draw[thick, blue] (c_1) -- (a_5);
    \draw[thick, blue] (c_1) -- (a_6);
    \draw[thick, blue] (w_i) -- (a_1);
    \draw[thick, blue] (a_1) -- (a_2);
    \draw[thick, green] (a_2) -- (a_3);
    \draw[thick, purple] (a_3) -- (a_4);
    \draw[thick, blue] (a_4) -- (a_5);
    \draw[thick, blue] (a_5) -- (a_6);
    \draw[thick, blue] (a_6) -- (u);
    \draw[thick, purple] (c_2) -- (a_4);
    \draw[thick, purple] (c_2) -- (a_3);
    \draw[thick, purple] (c_2) -- (b_1);
    \draw[thick, blue] (c_2) -- (b_2);
    \draw[thick, blue] (c_2) -- (b_3);
    \draw[thick, blue] (c_2) -- (b_4);
    \draw[thick, blue] (c_2) -- (b_5);
    \draw[thick, blue] (c_2) -- (b_6);
    \draw[thick, orange] (a_3) -- (b_1);
    \draw[thick, blue] (b_1) -- (b_2);
    \draw[thick, blue] (b_2) -- (b_3);
    \draw[thick, blue] (b_3) -- (b_4);
    \draw[thick, blue] (b_4) -- (b_5);
    \draw[thick, blue] (b_5) -- (b_6);
    \draw[thick, blue] (b_6) -- (a_4);   
\end{tikzpicture}
}

\newcommand{\HFour}{
\begin{tikzpicture}[scale=1.5]
    \node (v) at (-4,3) [circle, draw, fill=red!50, minimum size=10pt] {$v$};
    \node (u) at (1,0) [regular polygon, regular polygon sides=4, draw, fill=blue!50, minimum size=10pt] {$a_4$};
    \node (w_i) at (-1,0) [regular polygon, regular polygon sides=4, draw,fill=blue!50, minimum size=10pt] {$a_5$};
    \node (c_1) at (0,2) [regular polygon, regular polygon sides=4, shape border rotate=45, draw,fill=green!50, minimum size=5pt] {$2,c_a$};
    \node (a_1) at (-2,1) [regular polygon, regular polygon sides=4, draw,fill=blue!50, minimum size=5pt] {$a_6$};
    \node (a_2) at (-2,3) [circle, draw,fill=red!50, minimum size=5pt] {$u,a_7$};
    \node (a_3) at (-1,4) [regular polygon, regular polygon sides=4, shape border rotate=45, draw,fill=green!50, minimum size=5pt] {$1,a_0,b_0$};
    \node (a_4) at (1.2,4) [regular polygon, regular polygon sides=4, shape border rotate=45, draw,fill=green!50, minimum size=5pt] {$3,a_1,b_7$};
    \node (a_5) at (2,3) [regular polygon, regular polygon sides=4, draw,fill=blue!50, minimum size=5pt] {$a_2$};
    \node (a_6) at (2,1) [regular polygon, regular polygon sides=4, draw,fill=blue!50, minimum size=5pt] {$a_3$};
    \node (c_2) at (0,6) [regular polygon, regular polygon sides=4, shape border rotate=45, draw,fill=green!50, minimum size=5pt] {$4,c_b$};
    \node (b_1) at (-2.2,5.3) [regular polygon, regular polygon sides=4, shape border rotate=45, draw,fill=green!50, minimum size=5pt] {$5,\varsigma,b_1,\vartheta_{0}$};
    \node (b_2) at (-2,7.2) [regular polygon, regular polygon sides=6, draw,fill=purple!50, minimum size=5pt] {$b_2,\vartheta_{7}$};
    \node (b_3) at (-1,8) [regular polygon, regular polygon sides=4, draw, fill=blue!50, minimum size=5pt] {$b_3$};
    \node (b_4) at (1,8) [regular polygon, regular polygon sides=4, draw, fill=blue!50, minimum size=5pt] {$b_4$};
    \node (b_5) at (2,7) [regular polygon, regular polygon sides=4, draw,fill=blue!50, minimum size=5pt] {$b_5$};
    \node (b_6) at (2,5) [regular polygon, regular polygon sides=4, draw,fill=blue!50, minimum size=5pt] {$b_6$};
    \node (c_3) at (-4,6) [regular polygon, regular polygon sides=6, draw,fill=purple!50, minimum size=5pt] {$c_{\vartheta}$};
    \node (d_1) at (-3.4,4.1) [star, star points=5, draw, fill=orange!50, minimum size=5pt] {$\beta,\vartheta_{1}$};
    \node (d_2) at (-4.8,4.2) [regular polygon, regular polygon sides=4, draw,fill=blue!50, minimum size=5pt] {$\vartheta_{2}$};
    \node (d_3) at (-5.6,5.4) [regular polygon, regular polygon sides=4, draw,fill=blue!50, minimum size=5pt] {$\vartheta_{3}$};
    \node (d_4) at (-5.6,7) [regular polygon, regular polygon sides=4, draw,fill=blue!50, minimum size=5pt] {$\vartheta_{4}$};
    \node (d_5) at (-4.7,8) [regular polygon, regular polygon sides=4, draw,fill=blue!50, minimum size=5pt] {$\vartheta_{5}$};
    \node (d_6) at (-3,8) [regular polygon, regular polygon sides=4, draw,fill=blue!50, minimum size=5pt] {$\vartheta_{6}$};
    \draw[thick, blue] (w_i) -- (u);
    \draw[thick, green] (v) -- (a_3);
    \draw[thick, red] (a_2) -- (v);
    \draw[thick, blue] (c_1) -- (u);
    \draw[thick, blue](c_1) -- (w_i);
    \draw[thick, blue] (c_1) -- (a_1);
    \draw[thick, green] (c_1) -- (a_2);
    \draw[thick, green] (c_1) -- (a_3);
    \draw[thick, green] (c_1) -- (a_4);
    \draw[thick, blue] (c_1) -- (a_5);
    \draw[thick, blue] (c_1) -- (a_6);
    \draw[thick, blue] (w_i) -- (a_1);
    \draw[thick, blue] (a_1) -- (a_2);
    \draw[thick, green] (a_2) -- (a_3);
    \draw[thick, green] (a_3) -- (a_4);
    \draw[thick, blue] (a_4) -- (a_5);
    \draw[thick, blue] (a_5) -- (a_6);
    \draw[thick, blue] (a_6) -- (u);
    \draw[thick, green] (c_2) -- (a_4);
    \draw[thick, green] (c_2) -- (a_3);
    \draw[thick, green] (c_2) -- (b_1);
    \draw[thick, purple] (c_2) -- (b_2);
    \draw[thick, blue] (c_2) -- (b_3);
    \draw[thick, blue] (c_2) -- (b_4);
    \draw[thick, blue] (c_2) -- (b_5);
    \draw[thick, blue] (c_2) -- (b_6);
    \draw[thick, green] (a_3) -- (b_1);
    \draw[thick, purple] (b_1) -- (b_2);
    \draw[thick, blue] (b_2) -- (b_3);
    \draw[thick, blue] (b_3) -- (b_4);
    \draw[thick, blue] (b_4) -- (b_5);
    \draw[thick, blue] (b_5) -- (b_6);
    \draw[thick, blue] (b_6) -- (a_4);
    \draw[thick, purple] (c_3) -- (b_1);
    \draw[thick, purple] (c_3) -- (b_2);
    \draw[thick, purple] (c_3) -- (d_1);
    \draw[thick, blue] (c_3) -- (d_2);
    \draw[thick, blue] (c_3) -- (d_3);
    \draw[thick, blue] (c_3) -- (d_4);
    \draw[thick, blue] (c_3) -- (d_5);
    \draw[thick, blue] (c_3) -- (d_6);
    \draw[thick, orange] (b_1) -- (d_1);
    \draw[thick, blue] (d_1) -- (d_2);
    \draw[thick, blue] (d_2) -- (d_3);
    \draw[thick, blue] (d_3) -- (d_4);
    \draw[thick, blue] (d_4) -- (d_5);
    \draw[thick, blue] (d_5) -- (d_6);
    \draw[thick, blue] (d_6) -- (b_2);
    
\end{tikzpicture}
}

\newcommand{\HFive}{
\begin{tikzpicture}[scale=1.8]
    \node (v) at (-2.5,4.5) [circle, draw, fill=red!50, minimum size=10pt] {$v$};
    \node (u) at (1,0) [regular polygon, regular polygon sides=4, draw,fill=blue!50, minimum size=10pt] {$a_4$};
    \node (w_i) at (-1,0) [regular polygon, regular polygon sides=4, draw,fill=blue!50, minimum size=10pt] {$a_5$};
    \node (c_1) at (0,2) [regular polygon, regular polygon sides=4, shape border rotate=45, draw,fill=green!50, minimum size=5pt] {$2,c_a$};
    \node (a_1) at (-2,1) [regular polygon, regular polygon sides=4, draw,fill=blue!50, minimum size=5pt] {$a_6$};
    \node (a_2) at (-2,3) [circle, draw,fill=red!50, minimum size=5pt] {$u,a_7$};
    \node (a_3) at (-1.4,4) [regular polygon, regular polygon sides=4, shape border rotate=45, draw,fill=green!50, minimum size=5pt] {$1,a_0,b_0$};
    \node (a_4) at (1,4) [regular polygon, regular polygon sides=4, shape border rotate=45, draw,fill=green!50, minimum size=5pt] {$3,\varsigma,a_1,b_7,\vartheta_{7}$};
    \node (a_5) at (2,2.7) [regular polygon, regular polygon sides=4, draw,fill=blue!50, minimum size=5pt] {$a_2$};
    \node (a_6) at (2,1) [regular polygon, regular polygon sides=4, draw,fill=blue!50, minimum size=5pt] {$a_3$};
    \node (c_2) at (0,6) [regular polygon, regular polygon sides=4, shape border rotate=45, draw,fill=green!50, minimum size=5pt] {$4,c_b$};
    \node (b_1) at (-2,5.2) [regular polygon, regular polygon sides=4, draw,fill=blue!50, minimum size=5pt] {$b_1$};
    \node (b_2) at (-2,7) [regular polygon, regular polygon sides=4, draw,fill=blue!50, minimum size=5pt] {$b_2$};
    \node (b_3) at (-1,8) [regular polygon, regular polygon sides=4, draw,fill=blue!50, minimum size=5pt] {$b_3$};
    \node (b_4) at (1,8) [regular polygon, regular polygon sides=4, draw,fill=blue!50, minimum size=5pt] {$b_4$};
    \node (b_5) at (2,7) [regular polygon, regular polygon sides=4, draw,fill=blue!50, minimum size=5pt] {$b_5$};
    \node (b_6) at (2,5.3) [regular polygon, regular polygon sides=6, draw,fill=purple!50, minimum size=5pt] {$b_6,\vartheta_{0}$};
    \node (c_3) at (3.75,4.5) [regular polygon, regular polygon sides=6, draw,fill=purple!50, minimum size=5pt] {$c_{\vartheta}$};
    \node (d_1) at (3,6) [regular polygon, regular polygon sides=4, draw,fill=blue!50, minimum size=5pt] {$\vartheta_{1}$};
    \node (d_2) at (4.4,6) [regular polygon, regular polygon sides=4, draw,fill=blue!50, minimum size=5pt] {$\vartheta_{2}$};
    \node (d_3) at (5.2,5) [regular polygon, regular polygon sides=4, draw,fill=blue!50, minimum size=5pt] {$\vartheta_{3}$};
    \node (d_4) at (5.2,4) [regular polygon, regular polygon sides=4, draw,fill=blue!50, minimum size=5pt] {$\vartheta_{4}$};
    \node (d_5) at (4.4,3) [regular polygon, regular polygon sides=4, draw,fill=blue!50, minimum size=5pt] {$\vartheta_{5}$};
    \node (d_6) at (3,3) [star, star points=5, draw, fill=orange!50, minimum size=5pt] {$\beta, \vartheta_{6}$};
    \draw[thick, blue] (w_i) -- (u);
    \draw[thick, green] (v) -- (a_3);
    \draw[thick, red] (a_2) -- (v);
    \draw[thick, blue] (c_1) -- (u);
    \draw[thick, blue](c_1) -- (w_i);
    \draw[thick, blue] (c_1) -- (a_1);
    \draw[thick, green] (c_1) -- (a_2);
    \draw[thick, green] (c_1) -- (a_3);
    \draw[thick, green] (c_1) -- (a_4);
    \draw[thick, blue] (c_1) -- (a_5);
    \draw[thick, blue] (c_1) -- (a_6);
    \draw[thick, blue] (w_i) -- (a_1);
    \draw[thick, blue] (a_1) -- (a_2);
    \draw[thick, green] (a_2) -- (a_3);
    \draw[thick, green] (a_3) -- (a_4);
    \draw[thick, blue] (a_4) -- (a_5);
    \draw[thick, blue] (a_5) -- (a_6);
    \draw[thick, blue] (a_6) -- (u);
    \draw[thick, green] (c_2) -- (a_4);
    \draw[thick, green] (c_2) -- (a_3);
    \draw[thick, blue] (c_2) -- (b_1);
    \draw[thick, blue] (c_2) -- (b_2);
    \draw[thick, blue] (c_2) -- (b_3);
    \draw[thick, blue] (c_2) -- (b_4);
    \draw[thick, blue] (c_2) -- (b_5);
    \draw[thick, purple] (c_2) -- (b_6);
    \draw[thick, blue] (a_3) -- (b_1);
    \draw[thick, blue] (b_1) -- (b_2);
    \draw[thick, blue] (b_2) -- (b_3);
    \draw[thick, blue] (b_3) -- (b_4);
    \draw[thick, blue] (b_4) -- (b_5);
    \draw[thick, blue] (b_5) -- (b_6);
    \draw[thick, purple] (b_6) -- (a_4);
    \draw[thick, purple] (c_3) -- (a_4);
    \draw[thick, purple] (c_3) -- (b_6);
    \draw[thick, blue] (c_3) -- (d_1);
    \draw[thick, blue] (c_3) -- (d_2);
    \draw[thick, blue] (c_3) -- (d_3);
    \draw[thick, blue] (c_3) -- (d_4);
    \draw[thick, blue] (c_3) -- (d_5);
    \draw[thick, purple] (c_3) -- (d_6);
    \draw[thick, blue] (b_6) -- (d_1);
    \draw[thick, blue] (d_1) -- (d_2);
    \draw[thick, blue] (d_2) -- (d_3);
    \draw[thick, blue] (d_3) -- (d_4);
    \draw[thick, blue] (d_4) -- (d_5);
    \draw[thick, blue] (d_5) -- (d_6);
    \draw[thick, orange] (d_6) -- (a_4);
\end{tikzpicture}
}

\newcommand{\HSix}{
\begin{tikzpicture}[scale=1.2]
    \node (v) at (3,4) [circle, draw, fill=red!50, minimum size=10pt] {$v$};
    \node (u) at (1,0) [regular polygon, regular polygon sides=4, draw,fill=blue!50, minimum size=10pt] {$a_2$};
    \node (w_i) at (-1,0) [regular polygon, regular polygon sides=4, draw,fill=blue!50, minimum size=10pt] {$a_3$};
    \node (c_1) at (0,2) [regular polygon, regular polygon sides=6, draw,fill=purple!50, minimum size=5pt] {$c_a$};
    \node (a_1) at (-2,1) [regular polygon, regular polygon sides=4, draw,fill=blue!50, minimum size=5pt] {$a_4$};
    \node (a_2) at (-2,3) [regular polygon, regular polygon sides=4, draw,fill=blue!50, minimum size=5pt] {$a_5$};
    \node (a_3) at (-0.8,4) [star, star points=5, draw, fill=orange!50, minimum size=5pt] {$\beta, a_6$};
    \node (a_4) at (1,4) [circle, draw,fill=red!50, minimum size=5pt] {$u,a_7$};
    \node (a_5) at (2,3) [regular polygon, regular polygon sides=6, draw,fill=purple!50, minimum size=5pt] {$a_0$};
    \node (a_6) at (2,1) [regular polygon, regular polygon sides=4, draw,fill=blue!50, minimum size=5pt] {$a_1$};   
    \draw[thick, blue] (w_i) -- (u);
    \draw[thick, red] (v) -- (a_4);
    \draw[thick, purple] (v) -- (a_5);
    \draw[thick, blue] (c_1) -- (u);
    \draw[thick, blue] (c_1) -- (w_i);
    \draw[thick, blue] (c_1) -- (a_1);
    \draw[thick, blue] (c_1) -- (a_2);
    \draw[thick, purple]  (c_1) -- (a_3);
    \draw[thick, purple]  (c_1) -- (a_4);
    \draw[thick, purple]  (c_1) -- (a_5);
    \draw[thick, blue]  (c_1) -- (a_6);
    \draw[thick, blue] (w_i) -- (a_1);
    \draw[thick, blue] (a_1) -- (a_2);
    \draw[thick, blue] (a_2) -- (a_3);
    \draw[thick, orange] (a_3) -- (a_4);
    \draw[thick, purple]  (a_4) -- (a_5);
    \draw[thick, blue]  (a_5) -- (a_6);
    \draw[thick, blue] (a_6) -- (u);
    \end{tikzpicture}
}

\newif\ifpublic
\publictrue

\ifpublic

\newcommand{\ignore}[1]{}

\else

\fi

\newif\ifarXiv
\arXivtrue

\title{On fractional triangle decompositions of random graphs}

\author{Ghaura Mahabaduge}
\author{Michael Simkin}
\address{Department of Mathematics, Massachusetts Institute of Technology, Cambridge, MA 02139, USA}
\email{gmahabaduge3@gatech.edu,michael.simkin@mail.huji.ac.il}

\thanks{GM is supported by NSF grant CCF-2227876. MS is supported by NSF grant DMS-2349024.}

\begin{document}

\begin{abstract}
    We prove that with high probability $G(n,p)$ with $p \geq n^{-4/11 + o(1)}$ admits a fractional triangle decomposition (FTD), i.e., a nonnegative weighting of its triangles such that for each edge, the total weight of the triangles containing it equals one. This improves on the state of the art, due to Delcourt, Kelly, and Postle, that $p \geq n^{-1/3+o(1)}$ suffices.

    The proof is algorithmic: Given $G \sim G(n,p)$, we first construct an approximate FTD by taking a uniform weighting of the triangles. We then use specialized gadgets to iteratively shift weights and obtain successively better approximations of an FTD.
\end{abstract}

\maketitle

\section{Introduction}

Graph decompositions are a fundamental topic in combinatorics. For instance, Kirkman's theorem \cite{kirkman1847} that the complete graph $K_n$ admits a triangle decomposition (equivalently, an order-$n$ Steiner triple system exists) iff $n \equiv 1,3 \bmod 6$ is among the cornerstones of combinatorial design theory.

In general, one asks for conditions under which a graph $G$ can be decomposed into edge disjoint copies of a graph $H$. We do not expect a comprehensive characterization; even determining whether a given $G$ has a \textit{triangle} decomposition is NP-hard \cite{holyer1981np}. Nonetheless, the last decade has seen the introduction of powerful tools to prove the existence of decompositions under mild assumptions. These include \textit{randomized algebraic constructions} \cite{keevash2014existence,keevash2018counting}, \textit{iterative absorption} \cite{knox2015edge,kuhn2013hamilton,glock2023existence}, and \textit{refined absorption} \cite{delcourt2024refined}.

The simplest case, mentioned above, is when $H$ is a triangle. Although the triangle decomposition problem is NP-hard in general, we might hope to know when the random graph $G(n,p)$ admits such a decomposition w.h.p.\footnote{We say that a sequence of events parameterized by $n$ occurs \textit{with high probability (w.h.p.)} if the probabilities of their occurrence tend to $1$ as $n\to\infty$.} Upon reflection, however, we see that $G(n,p)$ almost never admits a triangle decomposition (unless $p$ is so small that there are typically no edges). Indeed, if a graph admits a triangle decomposition then it must satisfy the \textit{divisibility conditions} that the number of edges is a multiple of three and every degree must be even. While the former event occurs with reasonably high probability, the latter event almost never occurs (unless the graph is empty). Therefore, we modify the question: When does $G(n,p)$ admit a decomposition into (i) a matching of the odd-degree vertices (so that the degrees in the remaining graph are all even), (ii) at most one cycle of length four or five (so that the number of remaining edges is a multiple of three), and (iii) a collection of edge disjoint triangles? We say that a graph with this property admits a \textit{triangle$^*$ decomposition}.

While not its main focus, a straightforward consequence of Keevash's seminal work on the existence of designs \cite{keevash2018counting} is that there exists a constant $c>0$ such that for $p \geq n^{-c}$ w.h.p.\ $G(n,p)$ admits a triangle$^*$ decomposition. The threshold for triangle-decomposability (as well as more general clique decompositions of random graphs) was explicitly considered by Delcourt, Kelly, and Postle \cite{delcourt2024clique}. They applied the machinery of refined absorption to show that for every $\varepsilon>0$, w.h.p.\ $G(n,n^{-1/3+\varepsilon})$ can be decomposed into triangles and a set of at most $3n$ edges. The bottleneck here is the existence of certain gadgets in $G(n,p)$ that are used to modify an approximate triangle decomposition (covering all but a $o(1)$-fraction of edges) to an optimal one (covering all but $O(n)$ edges). These gadgets are somewhat dense, and their existence is unlikely if $p \ll n^{-1/3}$.

In this paper we ask when $G \sim G(n,p)$ admits a \textit{fractional triangle decomposition (FTD)}, i.e., a nonnegative weighting of the triangles in $G$ such that for every edge in $G$, the total weight of the triangles containing it equals one. FTDs differ from (integral) triangle decompositions in two key aspects. First, determining whether a given graph $G$ admits an FTD is a linear program, and so can be decided efficiently in both theory and practice. This allows us to perform numerical experiments. We share results of such experiments in \cref{sec:conclusion}. Second, the divisibility conditions mentioned above are no longer necessary. For instance, $K_n$ with $n \geq 3$ admits an FTD even though the number of edges may not be a multiple of three nor the degrees even.

The fractional triangle decomposability of $G(n,p)$ was considered in the aforementioned paper of Delcourt, Kelly, and Postle \cite{delcourt2024clique}. They proved that their lower bound on $p$ for which $G(n,p)$ is triangle$^*$ decomposable applies also to fractional triangle decomposability. That is, for every $\varepsilon>0$, w.h.p.\ $G(n,n^{-1/3+\varepsilon})$ admits an FTD. The proof is very similar to the integral case. The main difference is that the extra flexibility afforded by the fractional relaxation is used to ``fractionally absorb'' the small number of edges that are not covered by the approximate triangle packing.

Our main theorem is an improved lower bound on $p$ for which $G(n,p)$ admits an FTD w.h.p.

\begin{theorem}\label{thm:main}
    Let $\varepsilon > 0$ be fixed and let $G \sim G(n,p)$ with $p \geq n^{-4/11+\varepsilon}$. W.h.p.\ $G$ admits a fractional triangle decomposition.
\end{theorem}

What allows this improvement is that our techniques are geared specifically for the fractional case. However, as opposed to other work on fractional spanning structures in random graphs and hypergraphs (such as thresholds for fractional perfect matchings \cite{krivelevich1996perfect,devlin2017perfect} and $H$-factors \cite{haber2007fractional}), we do not use linear programming duality. Instead, we give a direct construction of a sequence of triangle weightings whose limit is an FTD. Our approach shares the combinatorial flavor of constructions due to Dross \cite{dross2016fractional}, Barber, K\"uhn, Lo, Montgomery, and Osthus \cite{barber2017fractional}, Montgomery \cite{montgomery2017fractional}, and Postle and Delcourt \cite{delcourt2021progress}, who studied the \textit{minimum degree threshold} for fractional triangle (or larger clique) decomposability. It is worth mentioning that Bowditch and Dukes \cite{bowditch2019fractional} applied linear algebraic methods to the minimum degree problem. Our proof is quite similar to that of \cite[Lemma 8.11]{kwan2023substructures} due to Kwan, Sah, Sawhney, and the second author. Our method generalizes easily, and can plausibly be applied to other fractional decomposition problems. We discuss this further in \cref{sec:conclusion}.

It is natural to ask what can be expected in terms of a lower bound. First, we note that admitting a triangle decomposition (whether fractional or integral) is not a monotone graph property. For instance, \triangleAndIso admits a triangle decomposition but \triangleWithTail does not admit an FTD. Hence, it is not obvious that a threshold for this property exists at all (even conditioning on nonemptiness of the sampled graph). Moreover, the non-monotonicity seems to preclude the use of the connection between fractional expectation thresholds and thresholds, proved by Frankston, Kahn, Narayanan, and Park \cite{frankston2021thresholds}, and which has been at the heart of many recent advances in the threshold theory of random (hyper)graphs, such as the threshold above which a random $3$-uniform hypergraph contains a spanning Steiner triple system \cite{jain2024optimal,keevash2022optimal}.

Nevertheless, we believe that there is a threshold for the property that $G(n,p)$ is non-empty and admits an FTD, as we now explain. A necessary condition for fractional triangle decomposability is that every edge is contained in at least one triangle. It is straightforward to prove that $p = \sqrt{ 3 \log n / (2n) }$ is a (sharp) threshold for this property. In the study of thresholds in random graphs it is often the case that a simple, local obstruction for a property is the only one. For instance, $G(n,p)$ is connected w.h.p.\ as soon as it does not have isolated vertices. With this in mind, we make the following conjecture.

\begin{conjecture}\label{con:main conjecture}
    For every $\varepsilon>0$ and $p \geq (1+\varepsilon) \sqrt{\frac{3 \log n}{2n}}$, w.h.p.\ $G(n,p)$ admits a fractional triangle decomposition.
\end{conjecture}
We remark that Yuster \cite{yuster2007combinatorial} conjectured that this is, up to a multiplicative constant, the threshold above which $G(n,p)$ admits a decomposition into triangles and at most $3n$ edges.

In \cref{sec:conclusion} we share numerical evidence strongly supporting \cref{con:main conjecture}, and state a stronger hitting time conjecture (\cref{con:hitting time}). We also discuss a generalization of our method that we view as a plausible path towards a proof.

The remainder of the paper is organized as follows. In \cref{sec:notation} we establish notation. In \cref{sec:overview} we give an overview of our method. In \cref{sec:concentration} we collect measure concentration inequalities that are useful in analyzing $G(n,p)$. In \cref{sec:combinatorial properties} we establish a sequence of combinatorial properties that hold w.h.p.\ in $G(n,n^{-4/11+o(1)})$. In \cref{sec:proof of main theorem} we prove that every graph that satisfies the combinatorial properties from \cref{sec:combinatorial properties} admits an FTD, thus completing the proof of \cref{thm:main}. Finally, in \cref{sec:conclusion}, we discuss numerical experiments and further research directions.

\section{Notation}\label{sec:notation}

\subsection{Graphs}
For a graph $G$ we write $V(G)$ and $E(G)$ for its vertex and edge sets, respectively. We write $v(G) = |V(G)|$ and $e(G) = |E(G)|$. We write $T(G)$ for its set of triangles. For a vertex $v \in V$ we write $T_G(v)$ for the set of triangles in $G$ containing $v$. If vertices $a,b,c$ span a triangle we sometimes abuse notation and write $abc$ for the triangle $\{a,b,c\}$.

Given a vertex $v \in V(G)$ we write $\deg_G(v)$ for its degree. If $u\in V(G)$ we write $\codeg_G(u,v)$ for the number of mutual neighbors of $u$ and $v$. In all of the above, we sometimes omit subscripts if the graph is clear from context.

For a function $f:T(G)\to\R$ and a vertex $v \in V(G)$ we write $f(v)=\sum_{T \in T(v)}f(T)$. Similarly, for an edge $e \in E(G)$ we write $f(e) = \sum_{e \subseteq T \in T(G)} f(T)$.

For a triangle $T \in T(G)$ we define the function $\one_T:T(G) \to \R$ as the indicator of the triangle $T$. 

\subsection{Rooted extensions}\label{sec:Rooted extensions}

At many points in the proof we will need to count the number of ways a partial embedding of a graph can be extended to a full embedding. We introduce the following notation.

\begin{definition}
    Let $H$ and $G$ be graphs, let $S \subseteq V(H)$ and let $\varphi : S \to V(G)$ be injective. We write $\chi(H,G,S,\varphi)$ for the set of injections $\psi: V(H) \to V(G)$ that satisfy
    \begin{itemize}
        \item $\psi|_S = \varphi$ (i.e., $\psi$ extends $\varphi$) and

        \item for every edge $uv \in E(H)$ with at least one of $u$ or $v$ not in $S$, $\psi(u)\psi(v) \in E(G)$ (i.e., $\psi$ preserves the edges of $H$, except possibly those induced by $S$).
    \end{itemize}
    We also write $X(H,G,S,\varphi) = |\chi(H,G,S,\varphi)|$.
\end{definition}

\section{Proof overview}\label{sec:overview}

In this section we describe the algorithm used to prove the existence of an FTD. Suppose that we have a graph $G=(V,E)$ (a typical sample from $G(n,n^{-4/11+\varepsilon})$). We wish to prove that $G$ supports an FTD. Our proof has three main phases: We first construct an approximate FTD by taking a uniform weight of the triangles, such that the total weight is the same as that of an FTD. We then modify the weights so that the total weight of the triangles containing each vertex is the same as that of an FTD. Finally, we apply a weight shifting operator that successively improves (i.e., makes closer to $1$) the total weight of each edge. We will prove that this sequence of weightings converges to an FTD.

We now describe the steps in more detail.

\subsection*{Step 1: The initial weighting}
Let $\varphi:T(G) \to \R$ be the uniform weighting of the triangles in $G$ given by
\[
\varphi \equiv \frac{e(G)}{3|T(G)|}.
\]
Note that the total weight of $\varphi$, which satisfies $\sum_{T \in T(G)} \varphi(T) = e(G)/3$, is the same as the total weight of all triangles in an FTD of $G$, since if $\rho:T(G)\to[0,1]$ is an FTD then
\[
\sum_{T \in T(G)} \rho(T) = \frac{1}{3} \sum_{e \in E} \sum_{e \subseteq T \in T(G)} \rho(T) = \frac{1}{3}e(G),
\]
where the second equality follows from the definition of an FTD.

\subsection*{Step 2: Correcting vertex weights}
In this step we modify $\varphi$ to obtain a weighting $\varphi_0$ in which all vertex weights are the same as they would be in an FTD. Note that if $\rho$ is an FTD, then for every vertex $v \in V$ we have 
\[
\sum_{v \in T \in T(G)} \rho(T) = \frac{1}{2} \sum_{v \in e \in E} \sum_{e \subseteq T \in T(G)} \rho(T) = \frac{1}{2}\deg_G(v),
\]
where the second equality follows from the definition of an FTD.

\begin{figure}
    \centering
    \bowtieFigure
    \caption{A $(u,v)$-bowtie. The shaded triangles are given weight $1$ while the unshaded triangles are given weight $-1$.}
    \label{fig:bowtie}
\end{figure}

For the function $\varphi$ above this condition might fail. We will modify it so that the condition is satisfied. We introduce the following notation: Given a function $\rho:T(G) \to \R$ and a vertex $v \in V$, let
\[
\delta(\rho,v) \coloneqq \rho(v) - \frac{1}{2}\deg_G(v)
\]
be the \textit{vertex discrepancy}. We also define $\delta(v) \coloneqq \delta(\varphi,v)$. We say that a function $\rho$ is \textit{vertex balanced} if $\delta(\rho,v)=0$ for all $v \in V$.

Our goal is to construct a vertex balanced weighting $\varphi_0 \approx \varphi$ (that also has the same total weight as $\varphi$). To do this we define a gadget that shifts weight between two vertices, while leaving the total weight of each other vertex (as well as the total weight of all triangles) unchanged. Given $u,v \in V$ a $(u,v)$-bowtie is a copy of the graph in \cref{fig:bowtie}. Observe that if the triangles $ua_1a_2$ and $cb_1b_2$ are given weight $1$, while the triangles $ca_1a_2$ and $vb_1b_2$ are given weight $-1$, then the total weight of $u$ is $1$, the total weight of $v$ is $-1$, and the total weight of each other vertex is $0$. Furthermore, the total weight of all triangles is $0$. Thus, we can use these gadgets to shift weight from vertices with excess weight to vertices with deficits. We give the remaining details in \cref{sec:vertex balancing}.

\subsection*{Step 3: Correcting edge deficiencies}
We now wish to modify $\varphi_0$ so that the total weight of each edge is $1$ (which is the defining constraint of an FTD). To this end, we introduce a gadget that shifts weight away from a given edge to other edges.

For $k \geq 3$, we denote by $W_k$ the wheel graph with vertex set $\{w_0,w_1,\ldots,w_{k-1},c\}$ and edge set $\{cw_i : 0 \leq i \leq k-1 \} \cup \{ w_iw_{(i+1) \bmod k} : 0 \leq i \leq k-1 \}$. Given an edge $uv \in E$, a \textit{$(u,v)$-octagonal pinwheel} is a copy of $W_8$ with $uv$ on its perimeter (see \cref{fig:octagon}). Observe that given a $(u,v)$-octagonal pinwheel, if we give its triangles weights $\pm 1$ alternating around the cycle, then the total weight of each vertex is $0$, as is the total weight of all triangles. Furthermore, the only edges whose weight changes are those on the perimeter, with the changes $1$ or $-1$ depending on the parity of their distance to $uv$.

\begin{figure}
\centering

\octagon
\caption{A $(u,v)$-octagonal pinwheel. The shaded triangles are given weight $1$, while the unshaded triangles are given weight $-1$.}\label{fig:octagon}
\end{figure}

Using these gadgets, we are able to shift each edge's discrepancy away from itself to the other edges of the pinwheel. In order that the discrepancy of any one edge is not simply absorbed into a small number of other edges, instead of taking a single pinwheel on each edge we take the average of all of them. In this way, each edge ``absorbs'' roughly the same fraction of each other edge's discrepancy, and these contributions mostly cancel each other out. This results in an overall reduction in discrepancy. By iterating this process, we obtain a sequence that converges to an FTD.

\subsection{What details are missing?}
Completing the sketch above to a proof primarily involves showing that copies of the gadgets above are distributed fairly evenly in the graph. For instance, one must show that for every (distinct) $(u,v)$ there are $\approx n^5p^{10}$ (labeled) $(u,v)$-bowties. Similarly, one must show that every $uv \in E$ is contained in $\approx n^7p^{15}$ (labeled) $(u,v)$-octagonal pinwheels. This is fairly straightforward to show, for example using the machinery of Kim--Vu polynomials (which we indeed employ heavily).

An astute reader might notice that the properties above (concentration of the number of bowties and pinwheels) occur also for some $p \ll n^{-4/11}$. Why, then, do we require this lower bound on $p$? The answer is that in Step 3, we must control the interactions between the collection of pinwheels on distinct edges (and sometimes even on sets of edges). For instance, for every pair of disjoint edges $u_1v_1,u_2v_2 \in E$, there are $\approx n^{12}p^{29}$ pairs of octagonal pinwheels over $(u_1,v_1)$ and $(u_2,v_2)$ that share the edge furthest from the base. In \cref{sec:combinatorial properties} we list several dozen such properties (with some parameterized, so that the printed list is not quite as long) that we prove hold w.h.p.\ in $G(n,p)$. These interactions can be quite subtle, and for some of the desired properties to hold we need to assume that $p \gg n^{-4/11}$.

\subsection{Why these particular gadgets?}
It is natural to ask whether the bowtie and octagonal pinwheel might be replaced by other gadgets, with the goal of either simplifying the proof or proving a stronger result.

In the case of the bowtie, the gadget indeed can be replaced by the simpler graph \diamondGraph (with the northwest and southeast vertices labeled $u$ and $v$, respectively). However, for every vertex pair $u,v$ to appear in a copy of this gadget we would require $p \gg n^{-2/5}$. While this is covered by the density $p \gg n^{-4/11}$ considered in \cref{thm:main}, we choose to use the bowtie since it can be used in the theoretically optimal range $p \gg n^{-1/2}$.

Regarding the octagonal pinwheel, the graphs $W_{2k}$ (i.e., wheels with even-length cycles) are natural candidates to replace it. It is entirely plausible that for every $\varepsilon>0$ there is some $k=k(\varepsilon)$ such that using $W_{2k}$ instead of $W_8$ would give an FTD of $G(n,p)$, for $p \geq n^{-1/2+\varepsilon}$. Unfortunately, the difficulty in analyzing the algorithm seems to increase with $k$. In choosing $W_8$, we have attempted to strike a balance between the amount of technical detail in the analysis and the improvement to the range in which $G(n,p)$ is known to admit an FTD. We remark that using $W_4$ or $W_6$ might actually be sufficient to prove \cref{thm:main} (certainly, they appear in abundance even when $p = n^{-4/11}$), but the difficulty of the analysis actually seems greater than with $W_8$.

\section{Measure concentration}\label{sec:concentration}

In this section we collect measure concentration inequalities that will be useful in the sequel. We use the following version of the Chernoff bound.

\begin{theorem}\label{thm:chernoff}
Let $X = \sum_{i=1}^{n} X_i$, where $\{X_i\}_{i=1}^n$ are independent Bernoulli random variables. Let $\mu = \mb E [X]$. For every $\delta \in (0,1)$ we have
\[
\Pr[|X - \mu| \geq \delta \mu] \leq 2e^{-\mu \delta^2 / 3}.
\]
\end{theorem}

We also use the concentration of polynomials of random variables, as developed by Kim and Vu~\cite{kim-vu-2000polynomial}. We use the setup in \cite[Chapter 7.8]{alon_spencer_probabilistic_method_3rd}, which we now present.

Let $\mc H = (V(\mc H), E(\mc H))$ be an $n$-vertex hypergraph with nonnegative edge weights $\{w_e\}_{e \in E(\mc H)}$. Let $\{t_i\}_{i \in V(\mc H)}$ be mutually independent Bernoulli random variables with $\mb E [t_i] = p_i$. Consider the polynomial  $Y = \sum_{e \in E(\mc H)} w_e \Pi_{i \in e} t_i$. Let $k$ be the degree of Y (which is at most the maximum size of a hyperedge of $\mc H$).

For $A \subseteq V(\mc H)$ with $|A| \leq k$ we write $E_A \coloneqq \sum_{A \subseteq e \in E(\mc H)} w_e \prod_{i \in e \setminus A}p_i$. Let
\[
E_i \coloneqq \max_{A \in \binom{V(\mc H)}{i}} E_A.
\]
Set $\mu = \mb E[Y]$, set $E'=\max_{1\leq i \leq k } E_i$, and let $E=\max\{E',\mu\}$.

For each $k \geq 0$, let  $a_k = 8^k(k!)^{1/2} $, $d_k = 2e^2$. The next theorem bounds the probability that $Y$ deviates from its expectation.

\begin{theorem}[{\cite[Theorem 7.8.1]{alon_spencer_probabilistic_method_3rd}}]\label{thm:kim-vu}
    In the setup above, for every $\lambda >0$ we have
    \[
    \Pr [|Y - \mu| > a_k(EE')^{1/2}\lambda^k] < d_k n^{k-1} \exp(-\lambda).
    \]
\end{theorem}

We mainly use \cref{thm:kim-vu} to count rooted extensions in random graphs. To this end, we prove the following corollary of \cref{thm:kim-vu}.

\begin{corollary}\label{cor:kim-vu vertex density}
    Let $H$ be a graph and $S \subseteq V(H)$. Suppose that for some $\alpha >0$ and every $W \subseteq V(H) \setminus S$ we have
    \[
    e(H[S \cup W]) - e(H[S]) \leq \alpha |W|.
    \]
    Let $\varepsilon \in (0,1)$ be fixed and let $G \sim G(n,p)$ with $p \geq n^{-1/\alpha+\varepsilon}$. W.h.p.\ for every injection $\varphi : S \to V(G)$ we have
    \begin{enumerate}
        \item\label{itm:homomorphism count}
        \[
        X(H,G,S,\varphi) = \left( 1 \pm n^{-\varepsilon/6} \right)n^{v(H)-|S|} p^{e(H)-e(H[S])};
        \]
        and

        \item\label{itm:fixed vertex homomorphism count} for every vertex $v \in V(G) \setminus \varphi(S)$ the number of extensions $\phi \in \chi(H,G,S,\varphi)$ such that $v \in \phi(V(H))$ is at most $n^{-\varepsilon/7} X(H,G,S,\varphi)$.
    \end{enumerate}
\end{corollary}

\begin{proof}
    Let $\varphi : S \to V(G)$ be injective.

    Write $s=|S|$ and $m=v(H)-s$. If $e(H) - e(H[S])=0$, then $X(H,G,S,\varphi)=(n-s)!/(n-s-m)!=(1+O(1/n))n^m$, and the proportion of extensions containing any fixed vertex outside $\varphi(S)$ is $m/(n-s)$. Thus both claims hold for large $n$, since $0<\varepsilon<1$. So assume that $e(H) - e(H[S]) > 0$. From this it follows that $S \subsetneq V(H)$. This also implies that $\alpha \geq 1/2$. Indeed, if there is an edge $uv \subseteq V(H) \setminus S$ then taking $W=\{u,v\}$ we see that $2\alpha \geq 1$. Otherwise, there is an edge $uv$ with $u \in S$ and $v \in V(H) \setminus S$. In this case, taking $W= \{v\}$, we see that $\alpha \geq 1$.

    We prove \cref{itm:homomorphism count} with an appeal to \cref{thm:kim-vu}. Let $Y=X(H,G,S,\varphi)$. We set $V(\mc H) = \binom{V(G)}{2}$ (i.e., potential edges in $G$). For $E(\mc H)$ we take those edge sets corresponding to a rooted extension of $\varphi$ into the complete graph on $V(G)$ (but without the edges spanned by $\varphi(S)$). We give each hyperedge the weight of the number of labeled rooted extensions producing that edge set; we note that these weights are all equal.

    We observe that the polynomial $Y$ has degree $k \coloneqq e(H) - e(H[S])$. We begin by calculating
    \[
    \mu = \mb E[Y] = p^{e(H) - e(H[S])} \frac{(n-|S|)!}{(n-v(H))!} = \left( 1 \pm O(1/n) \right)n^{v(H)-|S|} p^{e(H) - e(H[S])}.
    \]
    We observe that by assumption, taking $W = V(H) \setminus S$, we have $e(H) - e(H[S]) \leq \alpha (v(H) - |S|)$. Therefore $\mu = \Omega(n^{v(H) - |S| - \alpha(v(H) - |S|)(1/\alpha - \varepsilon)}) = \Omega (n^{\varepsilon \alpha}) = \Omega(n^{\varepsilon/2})$.
    
    We now bound $E_i$, for $1 \leq i \leq k$. Let $A \subseteq V(\mc H)$ be a set of $i$ edges. If $E_A >0$ then there exists at least one rooted extension $\psi$ that contains $A$. Note that in this case, no edge in $A$ has both endpoints in $\varphi(S)$ (otherwise $E_A$ would be zero). By assumption, $A$ is incident to at least $|A|/\alpha$ vertices in $V(G) \setminus \varphi(S)$. As a consequence,
    \[
    E_A \leq v(H)^{2|A|} p^{e(H) - e(H[S]) - |A|} n^{v(H) - |S| - |A|/\alpha} = O \left( \mu n^{(1/\alpha -\varepsilon)|A| - |A|/\alpha} \right) = O \left( \mu n^{-\varepsilon} \right).
    \]
    It follows that $E' = O (\mu n^{-\varepsilon})$ and that $E = \mu = \Omega(n^{\varepsilon/2})$. As a consequence, $\sqrt{EE'} = O(\mu n^{-\varepsilon/2})$. Taking $\lambda = c n^{\varepsilon/(3k)}$ for a sufficiently small fixed $c>0$, \cref{thm:kim-vu} implies that with probability $1 - \exp(- \Omega(n^{\varepsilon/{3k}}))$, we have $|Y-\mu| \leq \frac12 n^{-\varepsilon/6}\mu$. Hence, $Y = (1 \pm n^{-\varepsilon/6}) n^{v(H) - |S|} p^{e(H) - e(H[S])}$ for large $n$.
    
    The proof of \cref{itm:homomorphism count} now follows from a union bound over the (fewer than) $n^{|S|}$ choices of $\varphi$.

    To prove \cref{itm:fixed vertex homomorphism count}, let $v \in V(G) \setminus \varphi(S)$. Let $Y_v = X(H,G-v,S,\varphi)$ be the number of extensions of $\varphi$ that do not use $v$. Then the number of extensions $\phi \in \chi(H,G,S,\varphi)$ that use $v$ is $Y-Y_v$.

    Observe that $G-v \sim G(n-1,p)$. Thus, the calculations in the proof of \cref{itm:homomorphism count} apply with minor changes and we may use \cref{thm:kim-vu} to conclude that $Y_v = (1 \pm O(n^{-\varepsilon /6})) n^{v(H) - |S|} p^{e(H) - e(H[S])}$ with exponentially small failure probability; expressing the mean in terms of $n$ rather than $n-1$ incurs only a relative $O(1/n)$ error. Hence, $Y-Y_v = O(n^{-\varepsilon/6})X(H,G,S,\varphi) \leq n^{-\varepsilon/7}X(H,G,S,\varphi)$ for large $n$, with exponentially small failure probability. The claim now follows from a union bound over the (fewer than) $n^{|S|+1}$ choices of $\varphi$ and $v$.
\end{proof}

A useful way to bound $\alpha$ from \cref{cor:kim-vu vertex density} is with the notion of degeneracy.

\begin{definition}\label{def:extension degeneracy}
    Let $H$ be a graph and $S \subseteq V(H)$. We say that $(H,S)$ is \textit{$k$-degenerate} if there is an ordering $v_1,\ldots,v_N$ of $V(H) \setminus S$ such that for every $1 \leq i \leq N$ the vertex $v_i$ is adjacent to at most $k$ vertices in $S \cup \{v_1,\ldots,v_{i-1}\}$.
\end{definition}

\begin{claim}\label{clm:degeneracy}
    Let $H$ be a graph and $S \subseteq V(H)$, and suppose that $(H,S)$ is \textit{$k$-degenerate}. Then, for every $W \subseteq V(H) \setminus S$ we have $e(H[S \cup W]) - e(H[S]) \leq k|W|$.
\end{claim}

\begin{proof}
    Let $v_1,\ldots,v_N$ be an ordering of $V(H) \setminus S$ witnessing the $k$-degeneracy of $(H,S)$. Let $W \subseteq V(H) \setminus S$. For each $i$, let $H_i = H[S \cup \{v_1,\ldots,v_{i}\}]$. By assumption we have
    \[
    e(H[W \cup S]) - e(H[S]) \leq \sum_{v_i : v_i \in W} \deg_{H_i}(v_i) \leq k|W|,
    \]
    as claimed.
\end{proof}

\section{Combinatorial properties of $G(n,p)$}\label{sec:combinatorial properties}

\newcommand{\doubleLayerConcentratedSet}{{\mc P}}
\newcommand{\tripleLayerBaseSet}{{\mc T}}
\newcommand{\vertexBalanceBaseSet}{{\mc V}}
\newcommand{\tripleLayerSet}{{\mc Q}}

In this section we collect several properties that hold w.h.p.\ in $G(n,p)$, where $p \geq n^{-4/11+\varepsilon}$ and $\varepsilon>0$ is fixed. Most of the claims are to do with counting copies of rooted extensions. The motivation behind some properties may not be apparent before reading the proof of \cref{thm:main} in \cref{sec:proof of main theorem}. Thus, on first reading, the reader may wish to move on to \cref{sec:proof of main theorem} and refer to the current section as needed.

Throughout the section, let $\varepsilon > 0$ be fixed, let $p\geq n^{-4/11+\varepsilon}$, and let $G =(V,E) \sim G(n,p)$.

We begin with the standard claim that the degrees and codegrees in $G$ are concentrated around their respective expectations. The proof, which we omit, follows from straightforward applications of Chernoff's inequality and a union bound.

\begin{claim}\label{clm:degree and codegree concentration}
    The following hold w.h.p.
    \begin{enumerate}[{\bfseries{G\arabic{enumi}}}]
        \item\label{itm:deg concentration}\label{itm:first G} For every $v \in V$ we have $\deg_G(v) = (1 \pm n^{-0.2})np$; and

        \item\label{itm:codeg concentration} for all distinct $u,v \in V$ we have $\codeg_G(u,v) = (1 \pm n^{-0.1})np^2$.
        \setcounter{GProperties}{\value{enumi}}
    \end{enumerate}
\end{claim}

As previewed in \cref{sec:overview}, in \cref{sec:vertex balancing} we define an operation that takes an approximate FTD and shifts weight between vertices. To this end, we define the following gadget.

\begin{definition}\label{def:bowtie}
    Suppose that $u,v,a_1,a_2,b_1,b_2,c \in V$ are distinct and that
    \[
    ua_1, ua_2,a_1a_2,a_1c,a_2c, vb_1,vb_2, b_1b_2, b_1c, b_2c \in E
    \]
    (see \cref{fig:bowtie}).
    The function
    \[
    \psi_{u,v,a_1,a_2,b_1,b_2,c} = \one_{ua_1a_2} + \one_{cb_1b_2} - \one_{ca_1a_2} - \one_{vb_1b_2}
    \]
    is a \textit{$(u,v)$-bowtie}. We write $\mc B_{u,v}$ for the collection of all $(u,v)$-bowties in $G$.
\end{definition}

\begin{claim}\label{clm:bowtie counts}
    The following hold w.h.p.
    \begin{enumerate}[{\bfseries{G\arabic{enumi}}}]
        \setcounter{enumi}{\value{GProperties}}
        \item\label{itm:Bowtie Type concentration} For all distinct $u,v \in V$ we have $|\mc B_{u,v}| = (1 \pm n^{-0.02})n^5p^{10}\!$; and

        \item\label{itm:Triangles cant live in too many bowties} for every $T \in T(G)$ there are at most $48n^4p^7$ bowties $\psi$ on which $\psi(T) \neq 0$.
        \setcounter{GProperties}{\value{enumi}}
    \end{enumerate}
\end{claim}

\begin{proof}

We prove \ref{itm:Bowtie Type concentration} using \cref{cor:kim-vu vertex density} with $\alpha = 2$. Using notation from \cref{cor:kim-vu vertex density}, let $H$ be the bowtie graph in \cref{fig:bowtie} with $S = \{ u,v \}$. Let $W \subseteq V(H) \setminus S$. We wish to prove that $e(H[W \cup S]) - e(H[S]) \leq \alpha |W|$. Note that $S$ does not span any edges so that $e(H[W \cup S]) - e(H[S]) = e(H[W \cup S])$. We consider two cases.

\textit{Case 1:} If $c \notin W$ then $H[W \cup S]$ is contained in the union of the two triangles $ua_1a_2$ and $vb_1b_2$. It is easily seen that $e(H[W \cup S]) \leq \frac{3}{2}|W|$.

\textit{Case 2:} If $c \in W$ then let $W' \coloneqq W \setminus \{c\}$. we have $e(H[W \cup S]) = e(H[W' \cup S]) + |W'|$. By the first case we have $e(H[W' \cup S]) \leq \frac{3}{2}|W'|$. Therefore $e(H[W \cup S]) \leq \frac{3}{2}|W'| + |W'| = \frac{5}{2}(|W|-1)$. Since $|W| \leq 5$, the latter quantity is at most $2|W|$, as desired.

We now prove \ref{itm:Triangles cant live in too many bowties}, once again using \cref{cor:kim-vu vertex density} with $\alpha=2$. Let $H$ be the bowtie graph. Note that $H$ contains four triangles. Let $S$ be the vertices of one of these triangles. We will prove that $(H,S)$ is $2$-degenerate and apply \cref{clm:degeneracy}. We consider two cases.

\textit{Case 1: $S=\{u,a_1,a_2\}$.} In this case $c,b_1,b_2,v$ witnesses the $2$-degeneracy of $(H,S)$. The case $S=\{v,b_1,b_2\}$ is handled similarly.

\textit{Case 2: $S=\{c,a_1,a_2\}$.} In this case $u,b_1,b_2,v$ witnesses the $2$-degeneracy of $(H,S)$. The case $S=\{c,b_1,b_2\}$ is handled similarly.

\cref{cor:kim-vu vertex density} now implies that w.h.p.\ for every $T \in T(G)$ and every fixed mapping from $S$ to $T$ there are at most $2n^4p^7$ rooted extensions in $G$. Since there are four choices for $S$ and six ways to map $S$ to $T$, the claim follows.
\end{proof}

\begin{remark}
    Properties \cref{itm:deg concentration}--\cref{itm:Triangles cant live in too many bowties} hold even under the weaker assumption that $p \geq n^{-1/2+\beta}$, with $\beta>0$ a fixed constant, as long as one replaces the error terms with (say) $n^{-\beta/10}$. This is pertinent, since it means that \cref{lem:support vertex balanced} and its proof apply (with minor changes to error terms) in this density regime as well. This may be useful in an eventual proof of \cref{con:main conjecture}.
\end{remark}

In \cref{sec:edge adjustment} we define an ``edge adjustment'' operator that improves the edge discrepancies of approximate FTDs. The basic building block of this operator is the following gadget.

\begin{definition}\label{def:octagonal pinwheel}
    Recall that $W_8$ is the graph with vertex set $\{w_0,w_1,\ldots,w_7,c\}$ and edge set $\{cw_i : 0 \leq i \leq 7\} \cup \{ w_iw_{(i+1)\bmod8} : 0 \leq i \leq 7 \}$ (see \cref{fig:octagon}).

    Suppose that $u,v \in V$ are distinct. Let $\psi_{u,v}: \{ w_0,w_7 \} \to V$ map $w_0 \mapsto u$ and $w_7 \mapsto v$. A \textit{$(u,v)$-octagonal pinwheel} is an element of $\chi(W_8,G,\{w_0,w_7\},\psi_{u,v})$.

    Given a $(u,v)$-octagonal pinwheel $\psi$ with $uv \in E$, we define
    \[
    \Phi_\psi = \sum_{i=0}^7 (-1)^{i+1} \one_{\psi(c) \psi(w_i) \psi(w_{(i+1)\bmod8})}.
    \]

    Given an edge $e = uv \in E$ we write $\mc S(e) \coloneqq \chi(W_8,G,\{w_0,w_7\},\psi_{u,v}) \cup \chi(W_8,G,\{w_0,w_7\},\psi_{v,u})$ for the set of embeddings of $W_8$ that are either $(u,v)$-octagonal pinwheels or $(v,u)$-octagonal pinwheels.
\end{definition}

We begin by showing that the number of pinwheels over each edge is close to its expectation.

\begin{claim}\label{clm:pinwheel conc}
    W.h.p., the following hold
    \begin{enumerate}[{\bfseries{G\arabic{enumi}}}]
        \setcounter{enumi}{\value{GProperties}}
        \item\label{itm:pinwheel counts} For every $e \in E$ we have $|\mc S(e)| = (1 \pm n^{-\varepsilon/6}) 2n^7p^{15}$;

        \item\label{itm:pinwheel vertex avoidance} for every $e \in E$ and $x \in V \setminus e$ we have
        \[
        \left|\{\psi \in \mc S(e) : x \in \psi(V(W_8))\}\right| \leq n^{-\varepsilon/7}|\mc S(e)|;
        \]
        and

        \item\label{itm:pinwheel over triangle upper bound} for every triangle $T \in T(G)$, there are at most $50n^6p^{13}$ pinwheels $\psi \in \bigcup_{e \in E} \mc S(e)$ such that $\Phi_\psi(T) \neq 0$.
        \setcounter{GProperties}{\value{enumi}}
    \end{enumerate}
\end{claim}

\begin{proof}
    We first prove \cref{itm:pinwheel counts}. Recall that $|\mc S(\{u,v\})|$ counts octagonal pinwheels over $(u,v)$ and $(v,u)$; utilizing the reflective symmetry of $W_8$ the count is the same in both cases. Hence, \cref{itm:pinwheel counts} will follow from \cref{cor:kim-vu vertex density} with $H = W_8$, $S=\{w_0,w_7\}$ as long as we verify that the condition holds with $\alpha = 11/4$ (which it does with room to spare).

    Let $W \subseteq V(W_8) \setminus S$. We first observe that if $W=V(W_8) \setminus S$ then $e(W_8[W \cup S]) - e(W_8[S]) = e(W_8) - 1 = 15 \leq \frac{11}{4}\times 7 = \alpha |W|$.

    In all other cases $(W_8[W \cup S],S)$ is $2$-degenerate (see \cref{def:extension degeneracy}). Thus, \cref{clm:degeneracy} implies that $e(W_8[W \cup S]) - e(W_8[S]) \leq 2|W| \leq \alpha |W|$.

    Applying \cref{itm:fixed vertex homomorphism count} for both root orientations and summing the resulting bounds proves \cref{itm:pinwheel vertex avoidance}. A union bound then shows that for every $F \subseteq V \setminus e$, at most $|F|n^{-\varepsilon/7}|\mc S(e)|$ pinwheels in $\mc S(e)$ meet $F$.

    We now prove \cref{itm:pinwheel over triangle upper bound}. We begin by showing that the condition in \cref{cor:kim-vu vertex density} holds (with room to spare) when $\alpha = 11/4$, $H=W_8$, and $S$ is the vertex set of a triangle in $W_8$. W.l.o.g.\ we may assume that $S = \{c,w_0,w_1\}$. Let $W \subseteq V(W_8) \setminus S$. If $W \cup S = V(W_8)$ then $e(W_8[W \cup S]) - e(W_8[S]) = 16 - 3 = 13 \leq \frac{11}{4} \times 6 = \alpha |W|$.

    Otherwise $(W_8[W \cup S],S)$ is $2$-degenerate so \cref{clm:degeneracy} implies that $e(W_8[W \cup S]) - e(W_8[S]) \leq 2|W| \leq \alpha |W|$, as desired.

    It follows that w.h.p.\ for every injection $\varphi:S \to V(G)$ we have $X(W_8,G,S,\varphi) = (1 \pm n^{-\varepsilon/6}) n^6 p^{13}$. We deduce \cref{itm:pinwheel over triangle upper bound} by observing that there are eight triangles in $W_8$, and for each triangle $S$ in $W_8$ and every vertex triple $T \subseteq V(G)$, there are six ways to injectively map $S$ to $T$. Since each pinwheel $\psi$ with $\Phi_\psi(T) \neq 0$ corresponds to a rooted extension of a triangle to a copy of $W_8$, the number of such pinwheels is at most $(1+o(1)) 8 \times 6 \times n^6p^{13} < 50 n^6 p^{13}$, as claimed.
\end{proof}

In the next claim we show that for certain subgraphs of $W_8$ that can be found in $G$, the number of rooted extensions to a pinwheel is concentrated near its expectation.

\begin{definition}\label{def: Q_i}
    For every  $t \in [6]$ define  $Q_t = \{w_7,w_0,w_1,\ldots , w_{t-1},c\} \subseteq V(W_8)$.
\end{definition}

\begin{claim}
    W.h.p.\ the following hold.
    \begin{enumerate}[{\bfseries{G\arabic{enumi}}}]
        \setcounter{enumi}{\value{GProperties}}
        \item\label{itm:pinwheel completion} For every $t \in [5]$ and every injection  $\varphi:Q_t \to V(G)$  we have 
        \[
        X(W_8,G,Q_t,\varphi) = (1 \pm n^{-\eps/6})n^{7-t}p^{15-2t}.
        \]

        \item\label{itm:pinwheel 2-degeneracy} For every $t \in [6]$ and every injection $\varphi:\{w_0,w_7\} \to V(G)$  we have
        \[
        X(W_8 [Q_t],G,\{w_0,w_7\},\varphi) = (1 \pm n^{-\eps/6})n^{t}p^{2t}.
        \]
        \setcounter{GProperties}{\value{enumi}}
    \end{enumerate}
\end{claim}

\begin{proof}
We begin with \cref{itm:pinwheel completion}. Let $t \in [5]$. Note that $v(W_8)-v(Q_t)= 7-t$ and  $e(W_8)-e(W_8[Q_t])= 15-2t$. Therefore, it is sufficient to show that the condition for \cref{cor:kim-vu vertex density}  holds with $H=W_8$, $S=Q_t$, and $\alpha = 11/4$. Let $W \subseteq V(W_8) \setminus S$. If $W \cup S = V(W_8)$ then $e(W_8[W \cup S]) - e(W_8[S]) = 15 - 2t \leq \frac{11}{4} \times (7-t) = \alpha |W|$. Otherwise $(W_8[W \cup S],S)$ is $2$-degenerate so \cref{clm:degeneracy} implies that $e(W_8[W \cup S]) - e(W_8[S]) \leq 2|W| \leq \alpha |W|$, as desired.

We turn to \cref{itm:pinwheel 2-degeneracy}. Note that $|Q_t|-|\{w_0,w_7\}| = t$ and  $e(W_8[Q_t]) - e(W_8[\{w_0,w_7\}]) = 2t$. Therefore, it is sufficient to show that the condition in \cref{cor:kim-vu vertex density}  holds with $H=W_8[Q_t]$, $S=\{w_0,w_7\}$, and $\alpha = 11/4$. This follows from \cref{clm:degeneracy} and the fact that $(W_8[Q_t],S)$ is $2$-degenerate, as witnessed by the vertex ordering $c,w_1,w_2,\ldots,w_{t-1}$.
\end{proof}

As part of proving \cref{thm:main} we will need to understand ways in which pairs and triples of copies of $W_8$ interact. With this in mind, we define the following family of graphs.

\begin{definition}\label{def:Wij}
    For each $1 \leq i \leq 7$ we define the following graph $W(i)$. Take two copies of $W_8$, the first with vertices $a_0,\ldots,a_7$ in order around the cycle and central vertex $c_a$ and the second with vertices $b_0,\ldots,b_7$ around the cycle and central vertex $c_b$. Identify $a_{i}$ with $b_7$ and $a_{i-1}$ with $b_0$. Finally, add an additional vertex $d$ and the edges $da_0$ and $da_7$. (See \cref{fig:WThree} for an illustration of $W(3)$.)
\end{definition}

\begin{figure}
    \centering
    \wThree
    \caption{The graph $W(3)$.}
    \label{fig:WThree}
\end{figure}

In the proof of \cref{thm:main} we will need bounds on the number of rooted extensions of graphs in this family in which a perimeter edge as well as the vertices $d,a_7$ have been fixed. In order to keep notation compact, we define the following (disjoint) index sets, which form a partition of $[6] \times [7]$. The motivation behind the partition, as well as the set $\mc Q$ defined below, is tied to the methods we will use to analyze gadgets constructed from copies of $W_8$ in which certain vertices have been fixed. Gadgets parameterized by pairs in the same family will receive similar treatment.

\begin{definition}\label{def:TVP sets}

\begin{align*}
& \tripleLayerBaseSet = \{(1,2),(1,7)\},\\
& \vertexBalanceBaseSet = \{(1,1), (2,1), (6,7)\},\text{ and}\\
& \doubleLayerConcentratedSet = ([6] \times [7]) \setminus (\tripleLayerBaseSet \cup \vertexBalanceBaseSet).
\end{align*}
\end{definition}

\subsection*{Acknowledgment} The proof of \cref{clm:38 double-layer cases,clm:12 triple-layer cases,clm:color typicality} use python code to verify certain properties of graphs. We use the NetworkX library \cite{NetworkX} to construct and manage the graphs.

\begin{claim}\label{clm:38 double-layer cases}
    W.h.p.
    \begin{enumerate}[{\bfseries{G\arabic{enumi}}}]
        \setcounter{enumi}{\value{GProperties}}
        \item\label{itm:37 double-layer cases} for every $(i,j) \in \doubleLayerConcentratedSet$ and every injection $\varphi:\{d,a_7,b_{j-1},b_j\} \to V$ we have
        \[
        X(W(i),G,\{d,a_7,b_{j-1},b_j\},\varphi) = (1 \pm n^{-\varepsilon/6}) n^{13}p^{31}.
        \]
        \setcounter{GProperties}{\value{enumi}}
    \end{enumerate}
\end{claim}

\begin{remark}
    The lower bound $p = \omega(n^{-4/11})$ is necessary for \cref{clm:38 double-layer cases} to hold, and this is one of the sources of the lower bound on $p$ in \cref{thm:main}. To give just one example, consider the case $(i,j) = (2,7)$ (see \cref{fig:Case2.7}). The fixed vertices (shaded red in the figure) are $S = \{d,a_7,b_6,b_7\}$. For $W = \{c_a,a_0,a_1,c_b\}$ (the gray vertices in the figure) we have
    \[
    e(W(2)[W \cup S]) - e(W(2)[S]) = 11 = \frac{11}{4}|W|.
    \]
    Hence, for $p=O(n^{-4/11})$, we do not expect that every embedding of $S$ into $V$ will extend to a copy of $W(2)[W \cup S] \setminus W(2)[S]$ in $G(n,p)$.
\end{remark}

\begin{figure}[h]
    \centering
    \WTwoSeven
    \caption{The graph $W(2)$, corresponding to the case $(i,j) = (2,7)$. If the red vertices $\{ d,a_7,b_6,b_7 \}$ are fixed, the rooted extension to the gray vertices $\{ c_a,a_0,a_1,c_b \}$ has $11$ edges and four vertices, giving density $11/4$.}
    \label{fig:Case2.7}
\end{figure}

\begin{proof}[Proof of \cref{clm:38 double-layer cases}]
    It suffices to show that for $\alpha = 11/4$ the condition in \cref{cor:kim-vu vertex density} holds for every $(i,j) \in \doubleLayerConcentratedSet$. The python code attached to \ifarXiv{this arXiv submission }\else{the arXiv version of this paper \cite{mahabaduge2025fractional} }\fi verifies that this is the case: using brute force, the code verifies that for each $(i,j) \in \doubleLayerConcentratedSet$, writing $S_{i,j}\coloneqq \{d,a_7,b_{(j-1)\bmod8},b_j\}$, for every $W \subseteq V(W(i)) \setminus S_{i,j}$ we have $e[W(i)[W \cup S_{i,j}]] - e[W(i)[S_{i,j}]] \leq \frac{11}{4}|W|$.
\end{proof}

We turn our attention to $\tripleLayerBaseSet$. Here we consider concatenations of \textit{three} copies of $W_8$. We make the following definition.

\begin{definition}\label{def:W_ijm}
For $(i,j) \in  \tripleLayerBaseSet$ we define the graph $W(i,j)$ as follows. Take a copy of $W(i)$ and a copy of $W_8$, where we denote the central vertex of the (new) copy of $W_8$ as $c_{\vartheta}$ and the boundary vertices as $\vartheta_{0},\vartheta_{1}, \dots, \vartheta_{7}$ in order around the cycle. Identify $b_{j}$ with $\vartheta_7$ and $b_{j-1}$ with $\vartheta_0$. (See \cref{fig:Case1.2.6} for an illustration of $W(1,2)$.)
\end{definition}

\begin{figure}[h]
    \centering
    \WOneTwo
\caption{The graph $W(1,2)$ and the case $(1,2,6)$. The vertex set $\{d,a_7,\vartheta_5,\vartheta_{6}\}$ is highlighted in red.}\label{fig:Case1.2.6}
\end{figure}

\begin{definition}\label{def:triple layer base set def}
    We define
    \[
    \tripleLayerSet = \{(1,2,m) : 2 \leq m \leq 7\} \cup \{ (1,7,m) : m \in [6] \}.
    \]
\end{definition}
Observe that $\tripleLayerSet = \tripleLayerBaseSet \times [7] \setminus \{ (1,7,7), (1,2,1) \}$.

\begin{claim}\label{clm:12 triple-layer cases}
    W.h.p.\nopagebreak[4]
    \begin{enumerate}[{\bfseries{G\arabic{enumi}}}]
        \setcounter{enumi}{\value{GProperties}}
        \item\label{itm:12 triple-layer} for every $(i,j,m) \in \tripleLayerSet$ and every injection $\varphi:\{d,a_7,\vartheta_{m-1},\vartheta_m\} \to V$ we have
        \[
        X(W(i,j),G,\{d,a_7,\vartheta_{m-1},\vartheta_m\},\varphi) = (1 \pm n^{-\varepsilon/6}) n^{20}p^{46}.
        \]
        \setcounter{GProperties}{\value{enumi}}
    \end{enumerate}
\end{claim}

\begin{proof}

It suffices to show that for $\alpha = 11/4$ the condition in \cref{cor:kim-vu vertex density} holds for every $(i,j,m) \in \tripleLayerSet$. The python code attached to \ifarXiv{this arXiv submission }\else{the arXiv version of this paper \cite{mahabaduge2025fractional} }\fi verifies that this is the case.
\end{proof}

Recall that $\tripleLayerSet$ excludes the triples $(1,2,1),(1,7,7)$. Furthermore, we have not yet considered pairs or triples with $i=7$, nor the pairs in $\vertexBalanceBaseSet$. We will prove concentration results for the associated rooted extensions in a later claim. To set this up, we give the following definition.

\begin{definition}\label{def:color-typicality}
Let $H$ be a graph with vertex partition $V_R \sqcup V_G \sqcup V_O \sqcup V_{P} \sqcup V_{B}$ and let $q,\beta>0$.

An $n$-vertex graph $G' = (V',E')$ is \textit{$(H,V_R,V_G,V_O,V_{P},V_{B},q,\beta)$-typical} if the following hold:

\begin{enumerate}
    \item For every injection $\varphi: V_R \to V'$ we have
    \[
    X(H[V_G \cup V_R],G',V_R,\varphi) =(1\pm n^{-\beta})n^{|V_G|}q^{e(H[V_G \cup V_R])-e(H[V_R])}.
    \]

    \item For every injection $\varphi: V_R \cup V_G \cup V_O \to V'$ we have
    \[
    X(H,G', V_R \cup V_G \cup V_O,\varphi) =(1\pm n^{-\beta})n^{|V_P \cup V_B|}q^{e(H) - e(H[V_R \cup V_G \cup V_O])}.
    \]
\end{enumerate}
\end{definition}

The motivation behind this definition is that if typicality holds, then we have control over the number of copies of $H$ obtained in the multi-step process that begins with an arbitrary embedding of $V_R$, then extends to a copy of $H[V_G \cup V_R]$, then extends arbitrarily to an embedding of $H[V_R \cup V_G \cup V_O]$, and then extends to a complete copy of $H$.

We now define a collection of six graphs $H_1,\ldots,H_6$, each with an associated vertex partition as in \cref{def:color-typicality}. The graphs are defined in \cref{fig:H first}--\cref{fig:H last}. For each graph $H_i$, the sets $V_{R_i}$, $V_{G_i}$, $V_{O_i}$, $V_{P_i}$, and $V_{B_i}$ are, respectively, the sets of red (circular), green (diamond shaped), orange (star shaped), purple (hexagonal), and blue (square) vertices in the associated figure.

\begin{figure}
  \centering
    \HOne
    \caption{$H_1$: A copy of $W(2)$, used for the case $(i,j)=(2,1)$. Definition of the graph $H_1$, and the partition $V_R \sqcup V_G \sqcup V_O \sqcup V_P \sqcup V_B$ of its vertex set.}
    \label{fig:H first}\label{fig:HOne}
\end{figure}

\begin{figure}
    \centering
    \HTwo
    \caption{$H_2$: A copy of $W(6)$, used for the case $(i,j)=(6,7)$. Definition of the graph $H_2$, and the partition $V_R \sqcup V_G \sqcup V_O \sqcup V_P \sqcup V_B$ of its vertex set.}
    \label{fig:HTwo}
\end{figure}

\begin{figure}

    \centering
    \HThree
    \caption{$H_3$: A copy of $W(1)$, used for the case $(i,j)=(1,1)$. Definition of the graph $H_3$, and the partition $V_R \sqcup V_G \sqcup V_O \sqcup V_P \sqcup V_B$ of its vertex set.}
    \label{fig:HThree}
\end{figure}

\begin{figure}

    \centering
    \HFour
    \caption{$H_4$: A copy of $W(1,2)$, used for the case $(i,j,m)=(1,2,1)$. Definition of the graph $H_4$, and the partition $V_R \sqcup V_G \sqcup V_O \sqcup V_P \sqcup V_B$ of its vertex set.}
    \label{fig:HFour}
\end{figure}

\begin{figure}
    \centering
    \HFive
    \caption{$H_5$: A copy of $W(1,7)$, used for the case $(i,j,m)=(1,7,7)$. Definition of the graph $H_5$, and the partition $V_R \sqcup V_G \sqcup V_O \sqcup V_P \sqcup V_B$ of its vertex set.}
    \label{fig:HFive}

  \end{figure}

\begin{figure}
    \centering
    \HSix
    \caption{$H_6$, used for the single-wheel case $i=7$. Definition of the graph $H_6$, and the partition $V_R \sqcup V_G \sqcup V_O \sqcup V_P \sqcup V_B$ of its vertex set.}
    \label{fig:H last}\label{fig:HSix}
  \label{fig:grid6}
\end{figure}

\begin{claim}\label{clm:color typicality}
W.h.p.
\begin{enumerate}[{\bfseries{G\arabic{enumi}}}]
    \setcounter{enumi}{\value{GProperties}}
    \item\label{itm:color typicality}\label{itm:last G} for every $i \in [6]$, $G$ is \textit{$(H_i, V_{R_i}, V_{G_i}, V_{O_i}, V_{P_i}, V_{B_i}, p, \varepsilon/7)$-typical}.
    \setcounter{GProperties}{\value{enumi}}
\end{enumerate}
\end{claim}

\begin{proof}
    Let $i\in[6]$. To avoid excessive subscripting we set $V_R = V_{R_i}$, $V_G = V_{G_i}$, $V_O = V_{O_i}$, $V_P = V_{P_i}$, and $V_B = V_{B_i}$.
    
    We first prove the first condition of $(H_i,V_{R},V_{G},V_{O},V_{P},V_{B},p,\eps/7)$-typicality. To do so we apply \cref{cor:kim-vu vertex density} with $H = H_i[V_{G} \cup V_{R}]$, $S=V_{R}$, and $\alpha = 11/4$. For all $i$, $(H_i[V_{G} \cup V_{R}],V_{R})$ is $2$-degenerate. Indeed, in the figure associated to $H_i$, the numerical labels on the green vertices correspond to an ordering of $V_{G}$ that witnesses the $2$-degeneracy (note that this is vacuously true for $i=6$, since $V_{G_6}$ is empty). As a consequence, \cref{clm:degeneracy} and \cref{cor:kim-vu vertex density} imply that w.h.p.
    \[
    X(H_i[V_G \cup V_R],G,V_R,\varphi) = (1\pm n^{-\varepsilon/6})n^{|V_G|}p^{e(H[V_G \cup V_R])-e(H[V_R])}
    \]
    for every injection $\varphi : V_R \to V$.

    The second condition follows from \cref{cor:kim-vu vertex density}, with $H=H_i$, $S = V_R \cup V_G \cup V_O$, and $\alpha = 11/4$. The python code attached to \ifarXiv{this arXiv submission }\else{the arXiv version of this paper \cite{mahabaduge2025fractional} }\fi verifies that the condition holds.
\end{proof}

\section{Proof of \cref{thm:main}}\label{sec:proof of main theorem}

Let $\varepsilon > 0$ be fixed and let $p \geq n^{-4/11+\varepsilon}$. In the previous section we proved that w.h.p.\ $G(n,p)$ satisfies \cref{itm:first G}--\cref{itm:last G}. Thus, it suffices to prove that every graph that satisfies these conditions admits an FTD. Accordingly, for the remainder of this section we assume that $G$ is an $n$-vertex graph that satisfies \cref{itm:first G}--\cref{itm:last G} and that $n$ is large in terms of $\varepsilon$. We will prove that $G$ admits an FTD.

As outlined in \cref{sec:overview}, we describe an algorithm that constructs a sequence of weightings whose limit is an FTD. The algorithm has three stages. In the first stage we construct a constant-valued approximate FTD.

\subsection{Stage 1: Constructing an approximate FTD}\label{sec:uniform weighting}

We begin by defining a constant function that has the same total weight as an FTD of $G$ would have. Specifically, let $\varphi : T(G) \to \R$ be the constant function $\varphi \equiv |E|/(3|T(G)|)$. Observe that $\varphi$ has total weight $|E|/3$ --- the same total weight as any FTD of $G$.

\subsection{Stage 2: Correcting vertex defects}\label{sec:vertex balancing}

We now modify $\varphi$ so that every vertex has the same total weight as it would in an FTD.

\begin{definition}
    A function $\varphi_0:T(G)\to\R$ is \textit{vertex-balanced} if for every $v\in V(G)$ we have
    \[
    \sum_{T \in T_G(v)} \varphi_0(T) = \frac{1}{2} \deg_G(v).
    \]
\end{definition}

We define a vertex's \textit{defect} as the difference between its total weight under $\varphi$ and the weight it would have in an FTD. That is, for every $v \in V$ we define
\begin{equation}\label{eq:vertex defect def}
\delta(v) = \varphi(v) - \frac{1}{2}\deg_G(v).
\end{equation}
From the definitions it immediately follows that
\begin{equation}\label{eq:total vertex defect}
    \sum_{v \in V} \delta(v) = \sum_{v \in V} \left( |T_G(v)| \frac{|E|}{3|T(G)|} - \frac{1}{2} \deg_G(v) \right) = 0.
\end{equation}

We will use a weight-shifting ``gadget'' to adjust $\varphi$ in order to make all vertex defects zero.

\begin{lemma}\label{lem:support vertex balanced}
    $G$ supports a vertex-balanced function $\varphi_0$ that also satisfies: for every $T\in T(G)$ we have $\varphi_0(T) = (1 \pm n^{-\varepsilon/6})(np^2)^{-1}$.
\end{lemma}

In order to prove \cref{lem:support vertex balanced} we modify $\varphi$ from \cref{sec:uniform weighting} to obtain a vertex-balanced triangle weighting $\varphi_0$. As weight shifting gadgets, we use the bowties from \cref{def:bowtie}. Recall that if $\psi \in \mc B_{u,v}$ is a $(u,v)$-bowtie, then $\psi(u)=1$, $\psi(v)=-1$, and $\psi(w)=0$ for all vertices $w \neq u,v$.

\begin{proof}[Proof of \cref{lem:support vertex balanced}]
    Let $\varphi$ be the uniform weighting defined in \cref{sec:uniform weighting}. We use bowties to shift weight between the triangles to obtain a vertex balanced function. Define
    \[
    \varphi_0 \coloneqq \varphi - \frac{1}{n}\sum_{u \in V} \delta(u) \sum_{u \neq v \in V} \frac{1}{|\mc B_{u,v}|} \sum_{\psi \in \mc B_{u,v}} \psi.
    \]

    We will now show that $\varphi_0$ is vertex balanced. Let $w \in V$. Note that a bowtie $\psi \in \mc B_{u,v}$ gives non-zero weight to $w$ only if $w=u$ (in which case $\psi(w)=1$) or $w=v$ (in which case $\psi(w) = -1$). Therefore:
    
    \begin{align*}
    \varphi_0(w) & = \varphi(w) - \frac{1}{n}\delta(w) \sum_{w \neq v \in V} \frac{1}{|\mc B_{w,v}|} \sum_{\psi \in \mc B_{w,v}} \psi(w) - \frac{1}{n}\sum_{w \neq u \in V}\delta(u) \frac{1}{|\mc B_{u,w}|} \sum_{\psi \in \mc B_{u,w}} \psi(w)\\
    & = \varphi(w) - \left(\frac{1}{n}\delta(w) \sum_{w \neq v \in V} \frac{1}{|\mc B_{w,v}|} \sum_{\psi \in \mc B_{w,v}} 1 \right) - \left( \frac{1}{n}\sum_{w \neq u \in V}\delta(u) \frac{1}{|\mc B_{u,w}|} \sum_{\psi \in \mc B_{u,w}} (-1) \right)\\
    & = \varphi(w) - \left( \frac{1}{n}\delta(w) \sum_{w \neq v \in V}  1 \right) + \left(\frac{1}{n}\sum_{w \neq u \in V}\delta(u) \right)\\
    & = \varphi(w) - \frac{n-1}{n}\delta(w) + \frac{1}{n} \left( \sum_{u \in V} \delta(u) - \delta(w) \right)\\
    & = \varphi(w) - \delta(w) + \frac{1}{n} \sum_{u \in V}\delta(u)\\
    & \stackrel{\eqref{eq:total vertex defect}}{=} \varphi(w) - \delta(w)\\
    & \stackrel{\eqref{eq:vertex defect def}}{=} \frac{1}{2}\deg_G(w),
    \end{align*}
    implying that $\varphi_0$ is vertex-balanced.

    We will now prove that for every $T \in T(G)$ we have $\varphi_0(T) = (1 \pm n^{-\varepsilon/6})(np^2)^{-1}$. Let $T \in T(G)$. we have
    \begin{equation}\label{eq:vertex shift triangle inequality}
    |\varphi_0(T) - 1/(np^2)| \leq |\varphi(T) - 1/(np^2)| + |\varphi_0(T) - \varphi(T)|.
    \end{equation}
    We will bound the two terms on the RHS separately. Recall that $\varphi(T) = |E|/(3|T(G)|)$. By \ref{itm:deg concentration} we have $|E| = \frac{1}{2} \sum_{v \in V} \deg_G(v) = (1 \pm n^{-0.2}) \frac{1}{2} n^2p$. Additionally, \cref{itm:codeg concentration} implies that
    \[
    |T(G)| = \frac{1}{3} \sum_{uv \in E(G)} \codeg_G(u,v) = \frac{1}{3} |E| (1 \pm n^{-0.1})np^2.
    \]
    As a consequence
    \begin{equation}\label{eq:type 1 dev}
    \left| \varphi(T) - \frac{1}{np^2} \right|  = \left| \frac{|E|}{3|T(G)|} - \frac{1}{np^2} \right| \leq 2n^{-0.1} \frac{1}{np^2}.
    \end{equation}

    In order to bound $|\varphi(T) - \varphi_0(T)|$ recall that by definition
    \[
    \varphi(T) - \varphi_0(T) = \frac{1}{n} \sum_{u \in V}\delta(u) \sum_{u \neq v \in V} \frac{1}{|\mc B_{u,v}|} \sum_{\psi \in \mc B_{u,v}} \psi(T).
    \]
    By \cref{itm:Bowtie Type concentration} for every $u,v$ we have $|\mc B_{u,v}| \geq n^5p^{10}\!/2$ and by \cref{itm:Triangles cant live in too many bowties} there are at most $48 n^4 p^7$ bowties $\psi$ on which $\psi(T) \neq 0$ (and for those bowties we have $|\psi(T)|=1$). Hence
    \begin{equation}\label{eq:type 2 dev}
    |\varphi_0(T) - \varphi(T)| \leq \frac{1}{n} 
    \times 48 n^4p^7 \times \frac{2}{n^5p^{10}} \max \{|\delta(v)| : v \in V\} = \frac{96}{n^2p^3} \max \{|\delta(v)| : v \in V\}.
    \end{equation}
    To complete the argument we will bound $\max \{|\delta(v)| : v \in V\}$. Let $v \in V$. Recall that $|T_G(v)|$ is the number of triangles in $G$ containing $v$. we have
    \[
    |T_G(v)| = \frac{1}{2} \sum_{u : uv \in E} \codeg_G(u,v) \stackrel{\text{\cref{itm:deg concentration,itm:codeg concentration}}}{=} (1 \pm 2n^{-0.1}) \frac{1}{2} n^2p^3.
    \]
    Therefore (recalling that $\varphi \equiv |E|/(3|T(G)|) = (1 \pm 2n^{-0.1})(np^2)^{-1}$), we have that
    \[
    |\delta(v)| = \left| |T_G(v)| \frac{|E|}{3|T(G)|} - \frac{1}{2}\deg_G(v) \right| \stackrel{\text{\cref{itm:deg concentration}}}{=} \left| (1 \pm 4n^{-0.1}) \frac{1}{2}np - (1 \pm n^{-0.2}) \frac{1}{2}np \right| \leq 3n^{-0.1}np.
    \]
    Plugging this into \eqref{eq:type 2 dev} we see that
    \[
    |\varphi_0(T) - \varphi(T)| \leq \frac{96}{n^2p^3} \times 3 n^{-0.1}np \leq \frac{n^{-\varepsilon/5}}{np^2}.
    \]
    Combining the last inequality with \eqref{eq:type 1 dev} and \eqref{eq:vertex shift triangle inequality} we conclude that
    \[
    |\varphi_0(T) - 1/(np^2)| \leq \frac{n^{-\varepsilon/6}}{np^2},
    \]
    as desired.
\end{proof}

\subsection{Stage 3: Edge adjustment}\label{sec:edge adjustment}
Given a triangle weighting $\sigma : T(G) \to \R$, for each $e \in E(G)$ we define its \textit{edge discrepancy} as

\begin{equation}\label{eq:edge discrepancy}
\delta_e(\sigma) \coloneqq \sigma(e) - 1.
\end{equation}
This is the difference between the total weight of the edge under $\sigma$ and the weight of the edge in an FTD (which is always $1$). We also define
\[
\delta_\infty(\sigma) \coloneqq \max_{e\in E(G)} |\delta_e(\sigma)|.
\]

Observe that $\sigma$ is an FTD iff $\sigma \geq 0$ and $\delta_\infty(\sigma)=0$. In this section we define an ``edge adjustment'' operation that improves the edge discrepancies of approximate FTDs. We will show that repeatedly applying this operation to the weighting $\varphi_0$ constructed in the previous section yields a sequence of weightings that converge to an FTD.

The weight-shifting gadget we use is the octagonal pinwheel (see \cref{def:octagonal pinwheel}). Recall that for an edge $e = uv \in E(G)$, $\mc S(e)$ denotes the set of octagonal pinwheels $\psi$ over either $(u,v)$ or $(v,u)$. For each $e \in E(G)$ we define

\begin{equation}\label{eq:average pinwheel}
\Phi_e \coloneqq \frac{1}{|\mc S(e)|} \sum_{\psi \in \mc S(e)} \Phi_\psi.
\end{equation}

Observe that for every edge $e$, the total weight of triangles containing it under $\Phi_e$ is $1$ (i.e., $\Phi_e(e) = 1$). We further note that for every vertex $v$ and every octagonal pinwheel $\psi \in \mc S(e)$, we have $\Phi_\psi(v)=0$. In particular, if $\sigma$ is vertex-balanced then adding a linear combination of $\{\Phi_e\}_{e\in E}$ to $\sigma$ preserves the vertex balance.

Given a triangle weighting $\sigma: T(G) \to \R$, we define the triangle weighting 
\begin{equation}\label{eq:edge-adjustment}
F(\sigma)
\coloneqq
\sigma-\sum_{e\in E(G)}\delta_e(\sigma)\Phi_e.
\end{equation}

The idea behind this definition is that for each edge $e$, subtracting $\delta_e(\sigma)\Phi_e$ from $\sigma$ erases the discrepancy of $e$ and ``spreads it out'' among the other edges in the graph. Of course, we do not expect a single application of $F$ to result in an FTD, since each edge also ``absorbs'' some of the discrepancies from the other edges. However, as we will show, the contributions from the various edges mostly cancel each other out, resulting in an overall reduced discrepancy.

Let $\varphi_0$ be the vertex-balanced approximate FTD constructed in \cref{lem:support vertex balanced}. Define the sequence $\{\varphi_t\}_{t=0}^\infty$ by

\begin{equation}\label{eq:phi-recursion}
\varphi_t \coloneqq F(\varphi_{t-1}).
\end{equation}
Note that as observed above, this is a sequence of vertex-balanced weightings. In the next section we prove that the sequence converges, that the maximal edge-discrepancy of $\varphi_t$ decreases at a geometric rate, and that the functions $\varphi_t$ are all nonnegative. In particular, $\lim_{t\to\infty} \varphi_t$ is an FTD.

\subsection{Analyzing the sequence $\{\varphi_t\}_{t=0}^\infty$}\label{sec:reducing discrepancy}

\begin{claim}\label{clm:inductive master claim}
    The following hold for every $t$.
    \begin{enumerate}[{\bfseries{I\arabic{enumi}}}]
        \item \label{itm:weight dif of a triangle} $\norm{\varphi_t-\varphi_{t+1}}_\infty \leq \delta_\infty(\varphi_t) \frac{50}{np^2}$;

        \item \label{itm:max discrepancy at this stage in a triangle}  $\delta_\infty(\varphi_t) \leq 2^{-t}n^{-\varepsilon/8}$; and

        \item \label{itm:max discrepancy at this stage in a neigborhood}  for every $uv \in E(G)$ and the mutual neighborhood $z_1,\ldots,z_k$ of $u,v$ we have
        \[
        \left| \sum_{i=1}^k \delta_{u z_i}(\varphi_t) \right| \leq 2^{-t} np^2 n^{-\varepsilon/7}.
        \]
    \end{enumerate}
\end{claim}

Before proving \cref{clm:inductive master claim} we show how it implies that $G$ admits an FTD.

\begin{claim}
    The limit $\lim_{t\to\infty}\varphi_t$ exists and is a fractional triangle decomposition of $G$.
\end{claim}

\begin{proof}
    To see that $\lim_{t\to\infty}\varphi_t$ exists it suffices to observe that by \cref{itm:weight dif of a triangle,itm:max discrepancy at this stage in a triangle} the successive differences $\varphi_{t+1}-\varphi_{t}$ decrease geometrically in the $\ell_\infty$ norm.

    Let $\phi \coloneqq \lim_{t\to\infty}\varphi_t$. To prove that $\phi$ is an FTD we must prove that $\delta_\infty(\phi)=0$ and that $\phi$ is non-negative. The former fact follows from \cref{itm:max discrepancy at this stage in a triangle} since $\delta_\infty(\phi) = \lim_{t\to\infty}\delta_\infty(\varphi_t) = 0$ (noting that $\delta_\infty$ is continuous with respect to pointwise convergence). The latter fact follows from the triangle inequality, since for every triangle $T \in T(G)$ we have
    \[
    \phi(T) \geq \varphi_0(T) - \sum_{t=0}^\infty \norm{\varphi_t-\varphi_{t+1}}_\infty \stackrel{\text{\cref{itm:weight dif of a triangle,itm:max discrepancy at this stage in a triangle}}}{\geq} \varphi_0(T) - \sum_{t=0}^\infty \frac{50}{np^2}2^{-t}n^{-\varepsilon/8} = \varphi_0(T) - \frac{100}{np^2}n^{-\varepsilon/8}.
    \]
    By \cref{lem:support vertex balanced} we have $\varphi_0(T) \geq (np^2)^{-1}\!/2$. As a consequence
    \[
    \phi(T) \geq \frac{1}{2np^2} - \frac{100n^{-\varepsilon/8}}{np^2} \geq 0,
    \]
    as desired.
\end{proof}

We prove \cref{clm:inductive master claim} (mostly) inductively, as follows. We will first prove \cref{itm:weight dif of a triangle} directly. We will then show that if \cref{itm:max discrepancy at this stage in a triangle,itm:max discrepancy at this stage in a neigborhood} hold at time $t-1$ then \cref{itm:max discrepancy at this stage in a triangle} holds at time $t$. Then, we will directly show that \cref{itm:max discrepancy at this stage in a triangle} and \cref{itm:max discrepancy at this stage in a neigborhood} hold at times $0$, $1$, and $2$ (these are the base cases of the induction). Finally, we will prove that if \cref{itm:max discrepancy at this stage in a triangle} and \cref{itm:max discrepancy at this stage in a neigborhood} hold at times $t-3$, $t-2$, and $t-1$ then \cref{itm:max discrepancy at this stage in a neigborhood} holds at time $t$.

We now prove that \cref{itm:weight dif of a triangle} holds for all $t$.

\begin{proof}[Proof of \cref{itm:weight dif of a triangle}]

Let $t \geq 0$ and $T \in T(G)$. By \eqref{eq:edge-adjustment} and \eqref{eq:phi-recursion}, 
\[
|\varphi_t(T) - \varphi_{t+1}(T)|
= \left| \sum_{e \in E(G)} \delta_e(\varphi_t) \Phi_{e}(T) \right| \leq \delta_{\infty}(\varphi_t) \sum_{e \in E(G)} \left| \Phi_{e}(T) \right| \leq \delta_{\infty}(\varphi_t) \sum_{e \in E(G)} \sum_{\psi \in \mc S(e)} \frac{|\Phi_\psi(T)|}{|\mc S(e)|}.
\]
By \cref{itm:pinwheel counts} for every $e$ we have $|\mc S(e)| \geq n^7p^{15}$. Additionally, by \cref{itm:pinwheel over triangle upper bound} there are at most $50n^6p^{13}$ pinwheels $\psi$ such that $\Phi_\psi(T) \neq 0$ (and for each of these pinwheels we have $|\Phi_\psi(T)|=1$). Hence
\[
|\varphi_t(T) - \varphi_{t+1}(T)| \leq \delta_\infty(\varphi_t) \frac{50n^6p^{13}}{n^7p^{15}} = \frac{50}{np^2} \delta_\infty(\varphi_t),
\]
and the conclusion follows.
\end{proof}

We will now show that for $t\geq 1$, if \cref{itm:max discrepancy at this stage in a triangle} and \cref{itm:max discrepancy at this stage in a neigborhood} hold at time $t-1$ then \cref{itm:max discrepancy at this stage in a triangle} holds at time $t$. 

\begin{proof}[Proof that \cref{itm:max discrepancy at this stage in a triangle}$(t-1)$, \cref{itm:max discrepancy at this stage in a neigborhood}$(t-1)$ $\implies$ \cref{itm:max discrepancy at this stage in a triangle}$(t)$]
Let $t \geq 1$. We assume that $\delta_\infty(\varphi_{t-1}) \leq 2^{-(t-1)} n^{-\varepsilon/8}$ and that for every $uv \in E(G)$ and mutual neighborhood $z_1,\ldots,z_k$ of $u,v$ we have
\[
\left| \sum_{i=1}^k \delta_{uz_i}(\varphi_{t-1}) \right| \leq 2^{-(t-1)} np^2 n^{-\varepsilon/7}.
\]

Let $e= uv \in E(G)$. We will bound $\delta_{e}(\varphi_t)$. Recall that when we apply $F$ to $\varphi_{t-1}$, on the one hand edge $e$ loses weight $\delta_{e}(\varphi_{t-1})$, but on the other hand it ``absorbs'' some weight from the discrepancies of other edges. We will bound this incoming weight. By
\eqref{eq:edge discrepancy}, \eqref{eq:edge-adjustment}, and
\eqref{eq:phi-recursion}, we have
\begin{align*}
    \delta_e(\varphi_t) & = \varphi_t(e) - 1\\
    & = \delta_e(\varphi_{t-1}) + \varphi_t(e) - \varphi_{t-1}(e)\\
    & = \delta_e(\varphi_{t-1}) - \sum_{f \in E} \delta_f(\varphi_{t-1}) \Phi_f(e)\\
    & = \delta_e(\varphi_{t-1}) - \delta_e(\varphi_{t-1}) \Phi_e(e) - \sum_{e \neq f \in E} \delta_f(\varphi_{t-1}) \Phi_f(e).
\end{align*}
Recall that by construction, $\Phi_e(e) = 1$. Therefore
\begin{align*}
    \delta_e(\varphi_t) & = - \sum_{e \neq f \in E} \delta_f(\varphi_{t-1}) \Phi_f(e)\\
    & = - \sum_{e \neq f \in E} \frac{\delta_f(\varphi_{t-1})}{|\mc S(f)|} \sum_{\psi \in \mc S(f)} \Phi_\psi(e)\\
    & \stackrel{\text{\cref{itm:pinwheel counts}}}{=} - \sum_{e \neq f \in E} (1 \pm 2n^{-\varepsilon/6}) \frac{\delta_f(\varphi_{t-1})}{2n^7p^{15}} \sum_{\psi \in \mc S(f)} \Phi_\psi(e)\\
    & = - \frac{1}{2n^7p^{15}}  \sum_{e \neq f \in E} \delta_f(\varphi_{t-1}) \sum_{\psi \in \mc S(f)} \Phi_\psi(e) \pm 2n^{-\varepsilon/6} \frac{\delta_\infty(\varphi_{t-1})}{2n^7p^{15}} \sum_{e \neq f \in E} \sum_{\psi \in \mc S(f)} |\Phi_\psi(e)|,
\end{align*}
where the second equality follows from \eqref{eq:average pinwheel}.

We now consider the ways in which a pinwheel over $f$ can give non-zero weight to $e$. Recall that a pinwheel $\Phi_{\psi}$ gives non-zero weight to an edge $e'$ only if there is some $i$ such that $\{\psi(w_i),\psi(w_{i+1})\} = e'$ (where the indexing here is modulo $8$). Let $\tau: \{w_0,w_7\}\to V$ be the mapping $w_7 \mapsto u,w_0 \mapsto v$, and let $\sigma : \{w_0,w_7\} \to V$ be the mapping $w_7 \mapsto v,w_{0} \mapsto u$. Let $\chi \coloneqq \chi(W_8, G, \{w_0,w_{7}\}, \tau) \cup \chi(W_8, G, \{w_0,w_{7}\}, \sigma)$. For this bookkeeping we label each wheel so that $w_0w_7$ maps to the receiving edge $e$; the edge $\psi(w_{i-1})\psi(w_i)$ is the source whose discrepancy is transferred to $e$. Finally, define
\begin{equation}\label{eq:delta e i}
\delta(e,i) \coloneqq
\sum_{\psi \in \chi}
\delta_{\psi(w_i)\psi(w_{i-1})}(\varphi_{t-1}).
\end{equation}
It then follows that
\begin{equation}\label{eq:delta e vs delta e i}
\begin{split}
\delta_e(\varphi_t) & = \frac{1}{2n^7p^{15}} \sum_{i=1}^7 (-1)^{i+1} \delta(e,i) \pm 2n^{-\varepsilon/6} \frac{\delta_\infty(\varphi_{t-1})}{2n^7 p^{15}} \times 7 |\mc S(e)|\\
& \stackrel{\text{\cref{itm:pinwheel counts}}}{=} \frac{1}{2n^7p^{15}} \sum_{i=1}^7 (-1)^{i+1} \delta(e,i) \pm 2n^{-\varepsilon/6} \times 8 \delta_\infty(\varphi_{t-1}).
\end{split}
\end{equation}

We will now prove that $|\delta(e,i)|$ is small for every $i\in [4]$. This will suffice, since using the symmetries in $\sigma$, $\tau$, and $W_8$, we see that $\delta(e,7)=\delta(e,1),\delta(e,6)=\delta(e,2)$, and $\delta(e,5)=\delta(e,3)$.

\begin{figure}
    \centering
    \deltaIIllustration
    \caption{The embedding procedure corresponding to the bound on $\delta(e,3)$. Beginning with the embedding $\sigma$ mapping $w_0w_7$ to $uv$ (the red vertices), we extend $\sigma$ to an embedding $\psi_1$ of $W_8[Q_3]$ by embedding the green vertices. We then extend $\psi_1$ to an embedding of $W_8[Q_4]$ by embedding the orange vertex $w_3$ to a mutual neighbor of $\psi_1(c)$ and $\psi_1(w_2)$. Finally, we extend $\psi_2$ to an embedding $\psi$ of $W_8$ by embedding the blue vertices.}
    \label{fig:delta i is 3 illustration}
\end{figure}

Let $i \in [4]$. Recall the definition of the graphs $Q_i$ in \cref{def: Q_i}. We rewrite the sum in \eqref{eq:delta e i} by breaking each $\psi \in \chi$ into pieces: first we embed $W_8[Q_i]$, we then extend this embedding by a single vertex to $W_8[Q_{i+1}]$, and then extend this to a complete embedding of $W_8$ (see \cref{fig:delta i is 3 illustration} for an illustration of the case $i=3$). That is,
\begin{align*}
|\delta(e,i)| \leq & \left| \sum_{\psi_1 \in \chi(W_8[Q_i],G,\{w_0,w_7\},\tau)} \sum_{\psi_2 \in \chi(W_8[Q_{i+1}], G, Q_i,\psi_1)} \delta_{\psi_2(w_i)\psi_2(w_{i-1})} \sum_{\psi \in \chi (W_8, G, Q_{i+1}, \psi_2)} \Phi_{\psi}(e) \right|\\*
+ & \left| \sum_{\psi_1 \in \chi(W_8[Q_i],G,\{w_0,w_7\},\sigma)} \sum_{\psi_2 \in \chi(W_8[Q_{i+1}], G, Q_i,\psi_1)} \delta_{\psi_2(w_i)\psi_2(w_{i-1})} \sum_{\psi \in \chi (W_8, G, Q_{i+1}, \psi_2)} \Phi_{\psi}(e) \right|
\end{align*}
(the difference between the two lines is whether $\psi_1$ is an extension of $\tau$ or $\sigma$). By \cref{itm:pinwheel completion} we have $X(W_8, G, Q_{i+1}, \psi_2) = (1 \pm n^{-\varepsilon/6}) n^{6-i}p^{13-2i}$ for every choice of $\psi_2$ as above. Therefore
\begin{align*}
|\delta(e,i)| \leq & n^{6-i}p^{13-2i} \sum_{\psi_1 \in \chi(W_8[Q_i],G,\{w_0,w_7\},\tau)} \left| \sum_{\psi_2 \in \chi(W_8[Q_{i+1}], G, Q_i,\psi_1)} \delta_{\psi_2(w_i)\psi_2(w_{i-1})} \right|\\
+ & n^{6-i}p^{13-2i} \sum_{\psi_1 \in \chi(W_8[Q_i],G,\{w_0,w_7\},\sigma)} \left| \sum_{\psi_2 \in \chi(W_8[Q_{i+1}], G, Q_i,\psi_1)} \delta_{\psi_2(w_i)\psi_2(w_{i-1})} \right|\\
+ & n^{-\varepsilon/6} \times n^{6-i} p^{13-2i} \times \delta_\infty(\varphi_{t-1}) \sum_{\psi_1 \in \chi(W_8[Q_i],G,\{w_0,w_7\},\tau)} X (W_8[Q_{i+1}], G, Q_i,\psi_1)\\
+ & n^{-\varepsilon/6} \times n^{6-i} p^{13-2i} \times \delta_\infty(\varphi_{t-1}) \sum_{\psi_1 \in \chi(W_8[Q_i],G,\{w_0,w_7\},\sigma)} X (W_8[Q_{i+1}], G, Q_i,\psi_1) .
\end{align*}
We will now use the induction hypothesis to bound the inner sum. We first observe that for every $\psi_1$ as above, choosing the mapping $\psi_2$ is the same as choosing $\psi_2(w_{i})$, which is any mutual neighbor of $\psi_1(w_{i-1})$ and $\psi_1(c)$ that is not in the vertex set $\psi_1(w_7),\psi_1(w_0),\ldots,\psi_1(w_{i-2})$. Hence,
\[
\sum_{\psi_2 \in \chi(W_8[Q_{i+1}], G, Q_i,\psi_1)} \delta_{\psi_2(w_i)\psi_2(w_{i-1})}
\]
is the sum of edge-discrepancies of all but at most $i < 8$ mutual neighbors of $\psi_1(w_{i-1})$ and $\psi_1(c)$. By the induction hypothesis it follows that
\[
\left|\sum_{\psi_2 \in \chi(W_8[Q_{i+1}], G, Q_i,\psi_1)} \delta_{\psi_2(w_i)\psi_2(w_{i-1})}\right| \leq 2^{-(t-1)} np^2 n^{-\varepsilon/7} + 8\delta_\infty(\varphi_{t-1}) \leq 2^{2-t} np^2 n^{-\varepsilon/7}.
\]
Therefore, applying \cref{itm:codeg concentration},
\begin{align*}
|\delta(e,i)| \leq & n^{6-i}p^{13-2i} \times 2^{2-t} np^2 n^{-\varepsilon/7} \times X(W_8[Q_i],G,\{w_0,w_7\},\tau)\\
+ & n^{6-i}p^{13-2i} \times 2^{2-t} np^2 n^{-\varepsilon/7} \times X(W_8[Q_i],G,\{w_0,w_7\},\sigma)\\
+ & n^{-\varepsilon/6} \times n^{6-i} p^{13-2i} \times \delta_\infty(\varphi_{t-1}) \times X(W_8[Q_i],G,\{w_0,w_7\},\tau) \times 2np^2\\
+ & n^{-\varepsilon/6} \times n^{6-i} p^{13-2i} \times \delta_\infty(\varphi_{t-1}) \times X(W_8[Q_i],G,\{w_0,w_7\},\sigma) \times 2np^2.
\end{align*}
By \cref{itm:pinwheel 2-degeneracy} we have $X(W_8[Q_i],G,\{w_0,w_7\},\tau) = (1 \pm o(1)) n^{i}p^{2i}$, and the same with $\tau$ replaced by $\sigma$. Therefore
\[
|\delta(e,i)| \leq (1 \pm o(1)) n^{6-i}p^{13-2i} \times 2^{3-t} np^2 n^{-\varepsilon/7} \times n^{i}p^{2i} \leq 2^{4-t} n^7 p^{15} \times n^{-\varepsilon / 7}.
\]

By \eqref{eq:delta e vs delta e i} and the induction hypothesis it follows that
\[
|\delta_e(\varphi_t)| \leq (1 + o(1)) \frac{7}{n^7p^{15}} 2^{4-t} n^7p^{15} \times n^{-\varepsilon/7} \leq 2^{-t} n^{-\varepsilon/8},
\]
as desired.
\end{proof}

We next prove the base cases, which are that \cref{itm:max discrepancy at this stage in a triangle} and \cref{itm:max discrepancy at this stage in a neigborhood} hold for $t \leq 2$. 

\begin{proof}[Proof of the base cases]
We start by proving that \cref{itm:max discrepancy at this stage in a triangle} holds with $t=0$. Note that from \cref{lem:support vertex balanced} we have that for every $T \in T(G)$ we have $\varphi_0(T) = \frac{1}{np^2}(1 \pm n^{-\eps/6}) $. Let $e=uv \in E$. we have
\[
\varphi_0(e) = \codeg_G(u,v) (1 \pm n^{-\eps/6}) \frac{1}{np^2} \stackrel{\text{\cref{itm:codeg concentration}}}{=} (1 \pm n^{-0.1}) np^2 (1 \pm n^{-\eps/6}) \frac{1}{np^2} = 1 \pm n^{-\varepsilon/6.1}.
\]
Since this holds for every edge $e$ it follows that $\delta_\infty(\varphi_0) \leq n^{-\varepsilon/6.1}$. This is stronger than \cref{itm:max discrepancy at this stage in a triangle}, but we will use this stronger result later. 

Now we prove the $t=1$ case for \cref{itm:max discrepancy at this stage in a triangle}. By \cref{itm:weight dif of a triangle} we have $\norm{\varphi_1-\varphi_{0}}_\infty \leq \delta_{\infty}(\varphi_0) \frac{50}{np^2}$. As we have proved that $\delta_\infty(\varphi_0) \leq n^{-\varepsilon/6.1}$ it follows that $\norm{\varphi_0-\varphi_1}_\infty \leq \frac{50}{np^2}n^{-\varepsilon/6.1}$. Therefore, for every $e=uv \in E$ we have
\begin{align*}
|\delta_e(\varphi_1)| = \left| \sum_{e \subseteq T \in T(G)} \varphi_1(T) - 1 \right| & \leq |\delta_e(\varphi_0)| + \sum_{e \subseteq T \in T(G)} |\varphi_0(T) - \varphi_1(T)|\\
& \leq \delta_\infty(\varphi_0) + \codeg_G(u,v) \norm{\varphi_0-\varphi_1}_\infty\\
& \stackrel{\text{\cref{itm:codeg concentration}}}{\leq} n^{-\varepsilon/6.1} + 2np^2 \times \frac{50}{np^2} n^{-\varepsilon/6.1} \leq n^{-\varepsilon/6.2}.
\end{align*}
Since this holds for every $e$ it follows that $\delta_\infty(\varphi_1) \leq n^{-\varepsilon/6.2}$. In particular, this implies that \cref{itm:max discrepancy at this stage in a triangle} holds for $t=1$. From a similar argument we conclude that $\delta_\infty(\varphi_2) \leq n^{-\varepsilon/6.3}$. In particular \cref{itm:max discrepancy at this stage in a triangle} holds for $t \leq 2$.

Let $t \in \{0,1,2\}$. We proved above that $\delta_\infty(\varphi_t) \leq n^{-\varepsilon/6.3}$. We now prove that \cref{itm:max discrepancy at this stage in a neigborhood} holds at time $t$. Let $u,v \in V$ be distinct and let $z_1,\ldots,z_k$ be their mutual neighbors. we have
\[
\left| \sum_{i=1}^k \delta_{uz_i}(\varphi_t) \right| \leq \codeg_G(u,v) \delta_\infty(\varphi_t) \stackrel{\text{\cref{itm:codeg concentration}}}{\leq} 2np^2 \times n^{-\varepsilon/6.3} \leq 2^{-t} np^2 n^{-\varepsilon/7},
\]
as desired.
\end{proof}

Now we will prove that for $t\geq 3$ , if \cref{itm:max discrepancy at this stage in a triangle} and \cref{itm:max discrepancy at this stage in a neigborhood} hold at times $t-1$, $t-2$, and $t-3$ then \cref{itm:max discrepancy at this stage in a neigborhood} holds at time $t$.

\begin{proof}[Proof that \cref{itm:max discrepancy at this stage in a triangle}$(t-3,t-2,t-1)$, \cref{itm:max discrepancy at this stage in a neigborhood}$(t-3,t-2,t-1)$ $\implies$ \cref{itm:max discrepancy at this stage in a neigborhood}(t)] Let $u,v \in V$ be distinct and let $z_1,\ldots,z_k$ be their mutual neighbors. For $j \in [k]$ let $\chi_j = \chi(W_8,G,\{w_0,w_7\},\{w_7\mapsto u , w_0\mapsto z_j \})$. In words, $\chi_j$ is the set of embeddings of $W_8$ into $G$, where $w_7$ is mapped to $u$ and $w_0$ is mapped to $z_j$. For $i \in [7]$ we define

\begin{equation}\label{eq:delta u v i}
\delta(u,v;i) \coloneqq
\sum_{j=1}^{k} \sum_{\psi \in \chi_j}
\frac{\delta_{\psi(w_i) \psi(w_{i-1})}(\varphi_{t-1})}
{|\mc S(\psi(w_i) \psi(w_{i-1}))|}.
\end{equation}

Here again the receiving edge $uz_j$ is fixed as the root, and position $i$ specifies the source edge. Specifically,
\[
\left| \sum_{j=1}^k \delta_{uz_j}(\varphi_t) \right| = 2 \left| \sum_{i=1}^7 \delta(u,v;i)(-1)^{i+1} \right| \leq 2 \sum_{i=1}^7 |\delta(u,v;i)|
\]
(the factor $2$ is to account for the fact that $\chi_j$ fixes a root orientation of $uz_j$, but both orientations contribute to the discrepancy of $uz_j$). We will bound $|\delta(u,v;i)|$ for all $i \in [7]$. By
\eqref{eq:delta u v i} and \cref{itm:pinwheel counts}, it holds that
\begin{align*}
\delta(u,v;i) & \stackrel{\text{\cref{itm:pinwheel counts}}}{=} \sum_{j=1}^{k} \sum_{\psi \in \chi_j} \frac{ \delta_{\psi(w_i) \psi(w_{i-1})}(\varphi_{t-1})}{(1\pm n^{-\eps/6})2n^7p^{15}}\\
& = \frac{1}{2n^7p^{15}}\sum_{j=1}^{k} \sum_{\psi \in \chi_j} \delta_{\psi(w_i) \psi(w_{i-1})}(\varphi_{t-1}) \pm n^{-\varepsilon/6} \frac{\delta_\infty(\varphi_{t-1})}{n^7p^{15}} \sum_{j=1}^k |\chi_j|.
\end{align*}
Observe that for every $j$ we have $|\chi_j| = |\mc S(uz_j)|/2$. Therefore,
\[
n^{-\varepsilon/6} \frac{\delta_\infty(\varphi_{t-1})}{n^7p^{15}} \sum_{j=1}^k |\chi_j| \stackrel{\text{\cref{itm:pinwheel counts}}}{\leq} 2 n^{-\varepsilon/6} k \delta_\infty(\varphi_{t-1}) \stackrel{\text{\cref{itm:codeg concentration}}}{\leq} 3 n^{-\varepsilon/6} np^2 \delta_\infty(\varphi_{t-1}).
\]

Let
\[
\gamma(u,v,i) \coloneqq \frac{1}{2n^7p^{15}} \left| \sum_{j=1}^{k} \sum_{\psi \in \chi_j} \delta_{\psi(w_i) \psi(w_{i-1})}(\varphi_{t-1}) \right|.
\]
We will bound $\gamma(u,v,i)$ for each $i \in [7]$. We start from the case $i=7$. We will break down the sum in the definition of $\gamma(u,v,7)$ by the image of $w_6$ and $w_0$. To this end, for every $x\in N_G(u)$ and every $j \in [k]$ let $\chi_j(x)$ be the set of $\psi \in \chi_j$ such that $\psi(w_6) = x$. Recall that by definition $\psi(w_7)=u$ for every $\psi\in\chi_j$. We then have
\begin{align*}
\gamma(u,v;7) & = \frac{1}{2n^7p^{15}} \left| \sum_{x \in N_G(u)} \sum_{j=1}^{k} \sum_{\psi \in \chi_j(x)} \delta_{\psi(w_7) \psi(w_{6})}(\varphi_{t-1}) \right|\\
& = \frac{1}{2n^7p^{15}} \left| \sum_{x \in N_G(u)} \delta_{ux}(\varphi_{t-1}) \sum_{j=1}^{k} |\chi_j(x)| \right|.
\end{align*}
For $x \in N_G(u) \setminus\{v\}$ let $\chi_j^*(x) = \{\psi \in \chi_j(x) : v \notin \psi(V(W_8))\}$. Since $v \notin \{u,z_j\}$, \cref{itm:pinwheel vertex avoidance,itm:pinwheel counts,itm:codeg concentration} imply that the total absolute contribution from the excluded wheels with $x \neq v$ is at most
\begin{align*}
\frac{\delta_\infty(\varphi_{t-1})}{2n^7p^{15}} \sum_{j=1}^k \left|\{\psi \in \chi_j : v \in \psi(V(W_8))\}\right|
& \leq \frac{n^{-\varepsilon/7}\delta_\infty(\varphi_{t-1})}{2n^7p^{15}} \sum_{j=1}^k |\mc S(uz_j)|\\*
& \leq 2np^2n^{-\varepsilon/7}\delta_\infty(\varphi_{t-1}).
\end{align*}
We now note that for every $x \in N_G(u) \setminus\{v\}$, $\sum_{j=1}^k|\chi_j^*(x)|$ is the number of embeddings of the graph $H_6$ into $G$, where the vertices in $V_{R_6} \cup V_{G_6} \cup V_{O_6}$ have been fixed (recall that the graph $H_6$ and the sets $H_6,V_{R_6},V_{G_6},V_{O_6}$ are defined in \cref{fig:HSix}). By \cref{itm:color typicality} this is $(1 \pm n^{-\varepsilon/7})n^7p^{16}$. For $x=v$, $\sum_{j=1}^k|\chi_j(x)|$ is at most the number of extensions of $uv$ to a copy of $W_8$. By \cref{itm:pinwheel counts} this is at most $3n^7p^{15}$. Therefore
\begin{align*}
\gamma(u,v;7) & \leq p \left| \sum_{x \in N_G(u) \setminus \{v\}} \delta_{ux}(\varphi_{t-1}) \right| + n^{-\varepsilon/7}p \sum_{x \in N_G(u) \setminus \{v\}} \delta_{\infty}(\varphi_{t-1}) + 2\delta_\infty(\varphi_{t-1})\\
& \quad + 2np^2n^{-\varepsilon/7}\delta_\infty(\varphi_{t-1})\\
& \leq p \left| \sum_{x \in N_G(u)} \delta_{ux}(\varphi_{t-1}) \right| + p \delta_\infty(\varphi_{t-1}) + n^{-\varepsilon/7} p \delta_\infty(\varphi_{t-1}) \deg_G(u) + 2 \delta_\infty(\varphi_{t-1})\\
& \quad + 2np^2n^{-\varepsilon/7}\delta_\infty(\varphi_{t-1}).
\end{align*}
Since $\varphi_{t-1}$ is vertex-balanced the first sum is $0$. Additionally, by \cref{itm:deg concentration} we have $\deg_G(u) \leq 1.1 np$. Hence
\[
\gamma(u,v;7) \leq n^{-\varepsilon/7}4np^2\delta_\infty(\varphi_{t-1}).
\]

It follows that
\[
|\delta(u,v;7)| \leq \gamma(u,v;7) + 3n^{-\varepsilon/6} np^2 \delta_\infty(\varphi_{t-1}) \leq 5\delta_{\infty}(\varphi_{t-1}) np^2 n^{-\varepsilon/7}.
\]
By the induction hypothesis we have $\delta_\infty(\varphi_{t-1}) \leq 2^{-t+1}n^{-\varepsilon/8}$, so $|\delta(u,v;7)| \leq 10 \cdot 2^{-t}np^2n^{-15\varepsilon/56}$. Since $15/56>1/4$, for sufficiently large $n$ in terms of $\varepsilon$, uniformly in $t$, we have
\begin{equation}\label{eq:delta u v 7 bound}
|\delta(u,v;7)| \leq \frac{1}{100} 2^{-t} np^2 n^{-\varepsilon/4}.
\end{equation}

Note that in the proof above our strategy was to interpret $\delta(u,v;7)$ as a weighted sum of discrepancies of the form $\delta_{e}(\varphi_{t-1})$, i.e., we looked ``back in time'' a single step.

Now we will bound $|\delta(u,v;i)|$ for $i \in [6]$. For this our strategy is to write $\delta(u,v;i)$ as a weighted sum of discrepancies of the form $\delta_{e}(\varphi_{t-2})$ (i.e., look ``back in time'' two steps). The trick is to write $\delta(u,v;i)$ as a linear combination of $\{ \delta_{e}(\varphi_{t-1}) \}$ and then rewrite each term as a linear combination of $\{\delta_{e'}(\varphi_{t-2})\}$.

We make the following definition. Recall the graph $W(i)$ from \cref{def:Wij}. For every $\ell \in [k]$ let $\chi_{i,\ell} = \chi(W(i),G,\{d,a_0,a_7\},\{d\mapsto v, a_7\mapsto u,a_0 \mapsto z_\ell\})$. Let $j \in [7]$. We define

\begin{equation}\label{eq:delta u v i j}
\delta(u,v;i,j) \coloneqq
\sum_{\ell=1}^{k} \sum_{\psi \in \chi_{i,\ell}}
\frac{
\delta_{\psi(b_j) \psi(b_{j-1})}(\varphi_{t-2})
}{
|\mc S(\psi(b_j) \psi(b_{j-1}))|
|\mc S(\psi(a_i) \psi(a_{i-1}))|
}.
\end{equation}

Expanding the discrepancies in \eqref{eq:delta u v i} by \eqref{eq:edge-adjustment} gives pairs of individually injective wheels. The pairs that give an injective embedding of $W(i)$ contribute the terms in \eqref{eq:delta u v i j}. Thus, for $i \in [6]$,
\begin{equation}\label{eq:two layer overlap}
\delta(u,v;i) = 2 \sum_{j=1}^7 \delta(u,v;i,j)(-1)^{j+1} + R_i,
\end{equation}
where $R_i$ is the signed contribution of the remaining pairs. We will now bound $|R_i|$.

Write $B=n^7p^{15}$. A pair is excluded either because the first wheel contains $v$, or because the first wheel avoids $v$ and the second meets the first wheel or $v$ outside its prescribed root pair. By \cref{itm:pinwheel counts,itm:pinwheel vertex avoidance}, there are $O(B)$ choices for each wheel, while meeting a fixed forbidden set of at most eight vertices outside its root costs a factor $O(n^{-\varepsilon/7})$. Applying this bound at the first unwanted intersection gives $O(kn^{-\varepsilon/7}B^2)$ excluded pairs. It follows that
\begin{equation}\label{eq:two layer remainder}
|R_i| = O \left( k n^{-\varepsilon/7} B^2 \times \frac{\delta_\infty(\varphi_{t-2})}{B^2} \right) = O\left(2^{-t}np^2n^{-15\varepsilon/56}\right),
\end{equation}
where the last bound uses \cref{itm:codeg concentration} and the induction hypothesis.

By \eqref{eq:two layer overlap},
\[
|\delta(u,v;i)| \leq 2 \sum_{j=1}^7 |\delta(u,v;i,j)| + |R_i|.
\]
We will bound $\delta(u,v;i)$ by bounding each $|\delta(u,v;i,j)|$. To this end we define
\[
\gamma(u,v;i,j) \coloneqq \frac{1}{4n^{14}p^{30}} \left| \sum_{\ell=1}^{k} \sum_{\psi \in \chi_{i,\ell}} \delta_{\psi(b_j) \psi(b_{j-1})}(\varphi_{t-2}) \right|.
\]

We start by considering pairs in $\doubleLayerConcentratedSet$ (see \cref{def:TVP sets}). Let $(i,j) \in \doubleLayerConcentratedSet$. By \eqref{eq:delta u v i j}, it holds that
\begin{align*}
|\delta(u,v;i,j)| & \stackrel{\text{\cref{itm:pinwheel counts}}}{=} \left| \sum_{\ell=1}^{k} \sum_{\psi \in \chi_{i,\ell}} \frac{\delta_{\psi(b_j) \psi(b_{j-1})}(\varphi_{t-2})}{((1\pm n^{-\eps/6})2n^7p^{15})^2} \right|\\
& = \left| \sum_{\ell=1}^{k} \sum_{\psi \in \chi_{i,\ell}} \frac{\delta_{\psi(b_j) \psi(b_{j-1})}(\varphi_{t-2})}{(1\pm 3n^{-\eps/6})4n^{14}p^{30}} \right|\\
& \leq \frac{1}{4n^{14}p^{30}} \left| \sum_{\ell=1}^{k} \sum_{\psi \in \chi_{i,\ell}} \delta_{\psi(b_j) \psi(b_{j-1})}(\varphi_{t-2})   \right| + n^{-\eps/6}\frac{\delta_\infty (\varphi_{t-2})}{n^{14}p^{30}} \sum_{\ell=1}^k |\chi_{i,\ell}|\\
& = \gamma(u,v,i,j) + n^{-\eps/6}\frac{\delta_\infty (\varphi_{t-2})}{n^{14}p^{30}} \sum_{\ell=1}^k |\chi_{i,\ell}|.
\end{align*}

Note that we can bound $|\chi_{i,\ell}|$ by applying \cref{itm:pinwheel counts} twice, first on the edge $uz_\ell$ and then on $\psi(a_i)\psi(a_{i-1})$. It follows that $|\chi_{i,\ell}| \leq (3n^7p^{15})^2$ and therefore
\begin{equation}\label{eq:sum chi i ell bound}
\begin{split}
n^{-\eps/6}\frac{\delta_\infty (\varphi_{t-2})}{n^{14}p^{30}} \sum_{\ell=1}^k| \chi_{i,\ell} | & \leq  n^{-\eps/6}\frac{\delta_\infty (\varphi_{t-2})}{n^{14}p^{30}} \sum_{\ell=1}^k (3n^7p^{15})^2 = n^{-\eps/6}\frac{\delta_\infty (\varphi_{t-2})}{n^{14}p^{30}} (3n^7p^{15})^2 \codeg_G(u,v)\\
& \stackrel{\text{\cref{itm:codeg concentration}}}{\leq} 10n^{-\eps/6} np^2 \delta_\infty (\varphi_{t-2}).
\end{split}
\end{equation}

We will break down $\gamma(u,v;i,j)$ according to the image of $b_j,b_{j-1}$. That is,
\begin{align*}
\gamma(u,v;i,j) & = \frac{1}{4n^{14}p^{30}} \left| \sum_{\alpha\beta \in E} \sum_{\ell=1}^{k} \sum_{\substack{\psi \in \chi_{i,\ell}:\\ \psi(b_j)\psi(b_{j-1})=\alpha\beta}} \delta_{\alpha\beta}(\varphi_{t-2}) \right|\\
& = \frac{1}{4n^{14}p^{30}} \left| \sum_{\alpha\beta \in E} \delta_{\alpha\beta}(\varphi_{t-2}) \sum_{\ell=1}^{k} |\{ \psi \in \chi_{i,\ell} : \psi(b_j)\psi(b_{j-1}) = \alpha \beta \}| \right|.
\end{align*}
We now observe that as long as $\alpha,\beta \notin \{u,v\}$ we have
\begin{equation}\label{eq:W(8) count}
\begin{split}
\sum_{\ell=1}^{k} & |\{ \psi \in \chi_{i,\ell} : \psi(b_j)\psi(b_{j-1}) = \alpha \beta \}|\\
& = X(W(i),G,\{d, a_7,b_j,b_{j-1}\},\{ d \mapsto v, a_7 \mapsto u, b_j \mapsto \alpha, b_{j-1} \mapsto \beta \})\\
&\quad + X(W(i),G,\{d, a_7,b_j,b_{j-1}\},\{ d \mapsto v, a_7 \mapsto u, b_j \mapsto \beta, b_{j-1} \mapsto \alpha \})\\
& \stackrel{\text{\cref{itm:37 double-layer cases}}}{=} (1 \pm n^{-\varepsilon/6}) 2 n^{13}p^{31}.
\end{split}
\end{equation}
On the other hand, if $\{\alpha,\beta\} \cap \{u,v\} \neq \emptyset$ then $\sum_{\ell=1}^{k} |\{ \psi \in \chi_{i,\ell} : \psi(b_j)\psi(b_{j-1}) = \alpha \beta \}| = 0$, since only injective embeddings are counted. Therefore
\begin{align*}
\gamma(u,v;i,j) & \leq \frac{1}{4n^{14}p^{30}} \left( \left| \sum_{e \in E} \delta_e(\varphi_{t-2}) (1 \pm n^{-\varepsilon/6}) \right| + 2 (\deg_G(u) + \deg_G(v)) \delta_\infty(\varphi_{t-2}) \right) 2n^{13}p^{31}\\
& \stackrel{\text{\cref{itm:deg concentration}}}{\leq} \frac{p}{2n}\left| \sum_{e\in E}\delta_e(\varphi_{t-2}) \right| + n^{-\varepsilon/6} \frac{p}{2n} |E| \delta_\infty(\varphi_{t-2}) + \frac{p}{2n} 5np \delta_\infty(\varphi_{t-2})
\end{align*}
The degree terms appear since we applied \eqref{eq:W(8) count} to all edges $e \in E$ (which will allow us to use vertex balance), but then corrected for the edges incident to $u$ and $v$. Since the sum of discrepancies over all edges is zero (a consequence of the vertex-balance property), the first summand is zero. Hence
\[
\gamma(u,v;i,j) \leq \frac{p}{2n} \delta_\infty(\varphi_{t-2}) \left( n^{-\varepsilon/6} |E| + 5np \right) \stackrel{\text{\cref{itm:deg concentration}}}{\leq} \frac{p}{2n} \delta_\infty(\varphi_{t-2}) n^{-\varepsilon/6} n^2p \leq np^2 n^{-\varepsilon/6} \delta_\infty(\varphi_{t-2}).
\]

It follows that for every $(i,j) \in \doubleLayerConcentratedSet$ we have
\[
|\delta(u,v;i,j)| \stackrel{\eqref{eq:sum chi i ell bound}}{\leq} \gamma(u,v;i,j) + 10 np^2n^{-\eps/6} \delta_\infty(\varphi_{t-2}) \leq 11 np^2n^{-\eps/6} \delta_\infty(\varphi_{t-2}).
\]
By the induction hypothesis we have $\delta_\infty(\varphi_{t-2}) \leq 2^{-t+2}n^{-\varepsilon/8}$. Therefore for $(i,j) \in \doubleLayerConcentratedSet$
\[
|\delta(u,v;i,j)| \leq   \frac{1}{100} 2^{-t}np^2n^{-\eps/4}.
\]

We now consider $(i,j) \in \vertexBalanceBaseSet$ (see \cref{def:TVP sets}). As in the previous case we get
\[
|\delta(u,v;i,j)| \leq  \gamma(u,v;i,j) + 10 n^{-\eps/6} np^2 \delta_\infty (\varphi_{t-2}).
\] 

Before proceeding, we point out a connection between $\gamma(u,v;i,j)$ with $(i,j) \in \vertexBalanceBaseSet$ and the graphs $H_1,H_2,H_3$ shown in \cref{fig:HOne,fig:HTwo,fig:HThree}. We illustrate this with $H_2$ and $(i,j) = (6,7)$. In this case, by definition
\begin{align*}
\gamma(u,v;6,7) & = \frac{1}{4n^{14}p^{30}} \left| \sum_{\ell=1}^{k} \sum_{\psi \in \chi_{6,\ell}} \delta_{\psi(b_7) \psi(b_{7-1})}(\varphi_{t-2}) \right|\\
& \leq \frac{1}{4n^{14}p^{30}} \left| \sum_{\substack{xy \in E\\ x,y\notin\{u,v\}\\uy \in E}} \delta_{xy}(\varphi_{t-2}) X(H_2,G,\{ d,a_7,b_6,b_7 \},\{ d \mapsto v, a_7 \mapsto u, b_6 \mapsto x, b_7 \mapsto y \}) \right|\\
& + \frac{1}{4n^{14}p^{30}} \left| \sum_{\substack{xy \in E\\ x,y\notin\{u,v\}\\ux \in E}} \delta_{xy}(\varphi_{t-2}) X(H_2,G,\{ d,a_7,b_6,b_7 \},\{ d \mapsto v, a_7 \mapsto u, b_6 \mapsto y, b_7 \mapsto x \}) \right|.
\end{align*}
Similar equations can be written for $(i,j) = (2,1)$ and $H_1$, and $(i,j) = (1,1)$ and $H_3$. Thus, we will break down $\gamma(u,v;i,j)$ first according to the image of $V_R \cup V_G$ in the graphs $H_1,H_2,H_3$ shown in \cref{fig:HOne,fig:HTwo,fig:HThree} and then further break it down according to the image of $V_R \cup V_G \cup V_O$ in the same graphs.

Fix $(i,j) \in \vertexBalanceBaseSet$ and let $H = H_h$ with $h \in \{1,2,3\}$ according to the correspondence above. Let $V_G=V_{G_h},V_R=V_{R_h},V_O=V_{O_h},V_P=V_{P_h},V_B=V_{B_h}$ be as in the definition of $H_h$. Let
\[
A_{R,G} \coloneqq \chi(H[V_G \cup V_R],G,V_R, \{d\mapsto v, a_7\mapsto u  \}).
\]
For each $\rho \in A_{R,G}$ we define
\[
A_{R,G,O}(\rho) \coloneqq \chi(H[V_R \cup V_G \cup V_O], G, V_R \cup V_G, \rho)
\]
and for each $\rho' \in A_{R,G,O}(\rho)$ we define
\[
\chi(\rho') \coloneqq \chi(H,G,V_R \cup V_G \cup V_O, \rho').
\]
We can then write
\begin{align*}
    \gamma(u,v;i,j) & = \frac{1}{4n^{14}p^{30}} \left|  \sum_{\rho \in A_{R,G}} \sum_{\rho' \in A_{R,G,O}(\rho)} \sum_{\psi \in \chi(\rho') } \delta_{\psi(b_j) \psi(b_{j-1})}(\varphi_{t-2}) \right|.
\end{align*}
We now observe that by the definitions of $H$ and $W(i)$ the edge $b_jb_{j-1}$ is the same as the (orange) edge $\varsigma\beta$ (see \cref{fig:HOne,fig:HTwo,fig:HThree}). Hence, for every $\rho \in A_{R,G},\rho'\in A_{R,G,O}(\rho), \psi \in \chi(\rho')$ we have $\delta_{\psi(b_j)\psi(b_{j-1})}(\varphi_{t-2}) = \delta_{\rho'(\beta)\rho(\varsigma)}(\varphi_{t-2})$. As a consequence
\begin{align*}
    \gamma(u,v;i,j) & = \frac{1}{4n^{14}p^{30}} \left|  \sum_{\rho \in A_{R,G}} \sum_{\rho' \in A_{R,G,O}(\rho)} \delta_{\rho'(\beta)\rho(\varsigma)}(\varphi_{t-2}) |\chi(\rho')| \right|\\
    & \stackrel{\text{\cref{itm:color typicality}}}{=} \frac{1}{4n^{14}p^{30}} \left|  \sum_{\rho \in A_{R,G}} \sum_{\rho' \in A_{R,G,O}(\rho)} \delta_{\rho'(\beta)\rho(\varsigma)}(\varphi_{t-2}) (1 \pm n^{-\varepsilon/7}) n^{|V_P \cup V_B|}p^{e(H)-e(H[V_R \cup V_G \cup V_O])} \right|\\
    & \leq \frac{n^{|V_P \cup V_B|}p^{e(H)-e(H[V_R \cup V_G \cup V_O])}}{4n^{14}p^{30}} \left|  \sum_{\rho \in A_{R,G}} \sum_{\rho' \in A_{R,G,O}(\rho)} \left( \delta_{\rho'(\beta)\rho(\varsigma)}(\varphi_{t-2}) \pm n^{-\varepsilon/7}\delta_\infty(\varphi_{t-2}) \right) \right|\\
    & \leq \frac{n^{|V_P \cup V_B|}p^{e(H)-e(H[V_R \cup V_G \cup V_O])}}{4n^{14}p^{30}} \sum_{\rho \in A_{R,G}} \left| \sum_{\rho' \in A_{R,G,O}(\rho)} \delta_{\rho'(\beta)\rho(\varsigma)}(\varphi_{t-2}) \right|\\
    & \quad + \frac{n^{|V_P \cup V_B|}p^{e(H)-e(H[V_R \cup V_G \cup V_O])}}{4n^{14}p^{30}} |A_{R,G}| \max_{\rho \in A_{R,G}}|A_{R,G,O}(\rho)| n^{-\varepsilon/7} \delta_\infty(\varphi_{t-2})\\
    & \stackrel{\text{\cref{itm:color typicality}}}{\leq} \frac{n^{|V_P \cup V_B|}p^{e(H)-e(H[V_R \cup V_G \cup V_O])}}{4n^{14}p^{30}} \sum_{\rho \in A_{R,G}} \left| \sum_{\rho' \in A_{R,G,O}(\rho)} \delta_{\rho'(\beta)\rho(\varsigma)}(\varphi_{t-2}) \right|\\
    & \quad + np^2 n^{-\varepsilon/7} \delta_\infty(\varphi_{t-2}).
\end{align*}

We now turn our attention to the first term. Let $\rho \in A_{R,G}$. Observe that each $\rho' \in A_{R,G,O}(\rho)$ is determined by the image of $\beta$ (i.e., by the image of the orange vertex). Every vertex in $N_G(\rho(\varsigma)) \setminus \rho(V_R \cup V_G)$ can play this role. Thus
\[
\sum_{\rho' \in A_{R,G,O}(\rho)} \delta_{\rho'(\beta)\rho(\varsigma)}(\varphi_{t-2}) = \sum_{x \in N_G(\rho(\varsigma))} \delta_{\rho(\varsigma)x}(\varphi_{t-2}) \pm |V_R \cup V_G| \delta_\infty(\varphi_{t-2}).
\]
Since $\varphi_{t-2}$ is vertex balanced we have $\sum_{x \in N_G(\rho(\varsigma))} \delta_{\rho(\varsigma)x}(\varphi_{t-2}) = 0$. Therefore
\[
\left| \sum_{\rho' \in A_{R,G,O}(\rho)} \delta_{\rho'(\beta)\rho(\varsigma)}(\varphi_{t-2}) \right| \leq |V_R \cup V_G| \delta_\infty(\varphi_{t-2}) \leq 5  \delta_\infty(\varphi_{t-2}).
\]
As a consequence
\begin{align*}
& \frac{n^{|V_P \cup V_B|}p^{e(H)-e(H[V_R \cup V_G \cup V_O])}}{4n^{14}p^{30}} \sum_{\rho \in A_{R,G}} \left| \sum_{\rho' \in A_{R,G,O}(\rho)} \delta_{\rho'(\beta)\rho(\varsigma)}(\varphi_{t-2}) \right|\\
& \leq \frac{n^{|V_P \cup V_B|}p^{e(H)-e(H[V_R \cup V_G \cup V_O])}}{4n^{14}p^{30}} |A_{R,G}| \times 5  \delta_\infty(\varphi_{t-2}) \stackrel{\text{\cref{itm:color typicality}}}{\leq} 5 p \delta_\infty(\varphi_{t-2}) \leq np^2 n^{-\varepsilon/7} \delta_\infty(\varphi_{t-2}).
\end{align*}

It follows that for every $(i,j) \in \vertexBalanceBaseSet$ we have
\[
|\delta(u,v;i,j)| \leq \gamma(u,v;i,j) + 5np^2n^{-\eps/7} \delta_\infty(\varphi_{t-2}) \leq 7np^2n^{-\eps/7} \delta_\infty(\varphi_{t-2}).
\]

By the induction hypothesis we have $\delta_\infty(\varphi_{t-2}) \leq 2^{-t+2}n^{-\varepsilon/8}$. Therefore for $(i,j) \in \vertexBalanceBaseSet$
\[
|\delta(u,v;i,j)| \leq   \frac{1}{100} 2^{-t}np^2n^{-\eps/4}.
\]

Note that $\doubleLayerConcentratedSet \cup \vertexBalanceBaseSet$ covers all cases with $j \in [7]$ and $i \in \{2,3,4,5,6\}$. Therefore for all $i \in \{2,3,4,5,6\}$, we have
\begin{equation}\label{eq:delta most i}
|\delta(u,v;i)| \leq 2 \sum_{j=1}^7 |\delta(u,v;i,j)| + |R_i| \leq \frac{14}{100} 2^{-t}np^2n^{-\eps/4} + |R_i|.
\end{equation}

In light of \eqref{eq:delta u v 7 bound} it remains to bound $\delta(u,v,1)$. Note that $\doubleLayerConcentratedSet \cup \vertexBalanceBaseSet$ cover all pairs $(1,j)$ besides $j = 2$ and $j = 7$.

Just as in considering $(i,j)\in\doubleLayerConcentratedSet$ we looked at two ``layers'' of $W_8$, we now consider three layers. That is, our strategy is to write $\delta(u,v;i,j)$ as a weighted sum of discrepancies of the form $\delta_{e}(\varphi_{t-3})$.

We make the following definition. For $(i,j) \in \tripleLayerBaseSet$ recall the definition of the graph $W(i,j)$ in \cref{def:W_ijm}. For every $\ell \in [k]$ let $\chi_{i,j,\ell} = \chi(W(i,j),G,\{d,a_0,a_7\},\{d\mapsto v, a_7\mapsto u,a_0 \mapsto z_\ell\})$. We then define, for $m \in [7]$,

\begin{equation}\label{eq:delta u v i j m}
\delta(u,v;i,j,m) \coloneqq
\sum_{\ell=1}^{k}
\sum_{\psi \in \chi_{i,j,\ell}}
\frac{
\delta_{\psi(\vartheta_m)\psi(\vartheta_{m-1})}(\varphi_{t-3})
}{
|\mc S(\psi(\vartheta_m)\psi(\vartheta_{m-1}))|
|\mc S(\psi(b_j)\psi(b_{j-1}))|
|\mc S(\psi(a_i)\psi(a_{i-1}))|
}.
\end{equation}

Expanding the discrepancies in \eqref{eq:delta u v i j} once more gives, for $j \in \{2,7\}$,
\begin{equation}\label{eq:three layer overlap}
\delta(u,v;1,j) = 2 \sum_{m=1}^7 \delta(u,v;1,j,m)(-1)^{m+1} + R_{1,j},
\end{equation}
where $R_{1,j}$ is the signed contribution of triples that do not give an injective embedding of $W(1,j)$.

The first two wheels and $v$ already form an injective copy of $W(1)$, with $17$ vertices. The third wheel must avoid these vertices except for its two prescribed roots, so there are $15$ forbidden vertices, including $v$. Writing $B=n^7p^{15}$, by \cref{itm:pinwheel counts,itm:pinwheel vertex avoidance}, there are $O(kB^2)$ prefixes and $O(n^{-\varepsilon/7}B)$ excluded extensions for each prefix. Hence,
\begin{equation}\label{eq:three layer remainder}
|R_{1,j}| = O\left(kn^{-\varepsilon/7}B^3 \times \frac{\delta_\infty(\varphi_{t-3})}{B^3}\right)
= O\left(2^{-t}np^2n^{-15\varepsilon/56}\right).
\end{equation}
The last bound uses \cref{itm:codeg concentration} and the induction hypothesis.

By \eqref{eq:three layer overlap},
\[
|\delta(u,v;1,j)| \leq 2 \sum_{m=1}^7 |\delta(u,v;1,j,m)| + |R_{1,j}|.
\]
For $j \in \{2,7\}$ we will bound $|\delta(u,v;1,j)|$ by bounding each $|\delta(u,v;1,j,m)|$.

By \eqref{eq:delta u v i j m},
\begin{align*}
|\delta(u,v;i,j,m)| & = \left| \sum_{\ell=1}^{k} \sum_{\psi \in \chi_{i,j,\ell}} \frac{\delta_{\psi(\vartheta_m) \psi(\vartheta_{m-1})}(\varphi_{t-3})}{|\mc S(\psi(\vartheta_m) \psi(\vartheta_{m-1}))||\mc S(\psi(b_j) \psi(b_{j-1}))||\mc S(\psi(a_i) \psi(a_{i-1}))|}\right|\\
& \stackrel{\text{\cref{itm:pinwheel counts}}}{=} \left| \sum_{\ell=1}^{k} \sum_{\psi \in \chi_{i,j,\ell}} \frac{\delta_{\psi(\vartheta_m) \psi(\vartheta_{m-1})}(\varphi_{t-3})}{((1\pm n^{-\eps/6})2n^7p^{15})^3} \right|\\
& = \left| \sum_{\ell=1}^{k} \sum_{\psi \in \chi_{i,j,\ell}} \frac{\delta_{\psi(\vartheta_m) \psi(\vartheta_{m-1})}(\varphi_{t-3})}{(1\pm 4n^{-\eps/6})8n^{21}p^{45}} \right|\\
& \leq \frac{1}{8n^{21}p^{45}} \left| \sum_{\ell=1}^{k} \sum_{\psi \in \chi_{i,j,\ell}} \delta_{\psi(\vartheta_m) \psi(\vartheta_{m-1})}(\varphi_{t-3}) \right| + n^{-\eps/6}\frac{\delta_\infty (\varphi_{t-3})}{n^{21}p^{45}} \sum_{\ell=1}^k |\chi_{i,j,\ell}|.
\end{align*}

Note that we can bound $|\chi_{i,j,\ell}|$ by applying \cref{itm:pinwheel counts} thrice, first on the edge $u z_\ell$, then on $\psi(a_i) \psi(a_{i-1})$ and then on $\psi(b_j) \psi(b_{j-1})$. It follows that
\begin{align*}
n^{-\eps/6}\frac{\delta_\infty (\varphi_{t-3})}{n^{21}p^{45}} \sum_{\ell=1}^k| \chi_{i,j,\ell} | & \leq  n^{-\eps/6}\frac{\delta_\infty (\varphi_{t-3})}{n^{21}p^{45}} \sum_{\ell=1}^k (2n^7p^{15})^3 = n^{-\eps/6}\frac{\delta_\infty (\varphi_{t-3})}{n^{21}p^{45}} (2n^7p^{15})^3 \codeg_G(u,v)\\
& \stackrel{\text{\cref{itm:codeg concentration}}}{\leq} 9n^{-\eps/6} np^2 \delta_\infty (\varphi_{t-3}).
\end{align*}

We now define
\begin{equation}\label{eq:gamma def}
\gamma(u,v;i,j,m) \coloneqq \frac{1}{8n^{21}p^{45}} \left| \sum_{\ell=1}^{k} \sum_{\psi \in \chi_{i,j,\ell}} \delta_{\psi(\vartheta_m) \psi(\vartheta_{m-1})}(\varphi_{t-3})   \right|.
\end{equation}
We will bound $\gamma(u,v;i,j,m)$ for $(i,j,m) \in \{(1,2),(1,7)\} \times [7]$.

We start with the case $(i,j,m) \in \{ (1,2,1) , (1,7,7) \}$. The proof is similar to the proof for $\vertexBalanceBaseSet$. Specifically, we will exploit the connection between $\gamma(u,v;1,2,1)$ and embeddings of $H_4$ in $G$ and between $\gamma(u,v;1,7,7)$ and embeddings of $H_5$ (recalling the definitions in \cref{fig:HFour,fig:HFive}).

Let $(i,j,m) \in \{ (1,2,1) , (1,7,7) \}$ and let $H \in \{H_4,H_5\}$ be the associated graph, together with the partition of its vertex set into $V_R,V_G,V_O,V_B,V_P$. We will break down $\gamma(u,v;i,j,m)$ first according to the image of $V_R \cup V_G$ and then according to the image of $V_R \cup V_G \cup V_O$.

Let
\[
A_{R,G} \coloneqq \chi(H[V_G \cup V_R],G,V_R, \{d\mapsto v, a_7\mapsto u \}).
\]
For each $\rho \in A_{R,G}$ we define
\[
A_{R,G,O}(\rho) \coloneqq \chi(H[V_R \cup V_G \cup V_O], G, V_R \cup V_G, \rho)
\]
and for each $\rho' \in A_{R,G,O}(\rho)$ we define
\[
\chi(\rho') \coloneqq \chi(H,G,V_R \cup V_G \cup V_O, \rho').
\]
We can then write
\begin{align*}
    \gamma(u,v;i,j,m) & = \frac{1}{8n^{21}p^{45}} \left|  \sum_{\rho \in A_{R,G}} \sum_{\rho' \in A_{R,G,O}(\rho)} \sum_{\psi \in \chi(\rho') } \delta_{\psi(\vartheta_m) \psi(\vartheta_{m-1})}(\varphi_{t-3}) \right|.
\end{align*}
We now observe that by the definitions of $H$ and $W(i,j)$ the edge $\vartheta_m\vartheta_{m-1}$ is the same as the (orange) edge $\varsigma\beta$ (see \cref{fig:HFour,fig:HFive}). Hence, for every $\rho \in A_{R,G},\rho'\in A_{R,G,O}(\rho), \psi \in \chi(\rho')$ we have $\delta_{\psi(\vartheta_m)\psi(\vartheta_{m-1})}(\varphi_{t-3}) = \delta_{\rho'(\beta)\rho(\varsigma)}(\varphi_{t-3})$. As a consequence
\begin{align*}
    \gamma(u,v;i,j,m) & = \frac{1}{8n^{21}p^{45}} \left|  \sum_{\rho \in A_{R,G}} \sum_{\rho' \in A_{R,G,O}(\rho)} \delta_{\rho'(\beta)\rho(\varsigma)}(\varphi_{t-3}) |\chi(\rho')| \right|\\
    & \stackrel{\text{\cref{itm:color typicality}}}{=} \frac{1}{8n^{21}p^{45}} \left|  \sum_{\rho \in A_{R,G}} \sum_{\rho' \in A_{R,G,O}(\rho)} \delta_{\rho'(\beta)\rho(\varsigma)}(\varphi_{t-3}) (1 \pm n^{-\varepsilon/7}) n^{|V_P \cup V_B|}p^{e(H)-e(H[V_R \cup V_G \cup V_O])} \right|\\
    & \leq \frac{n^{|V_P \cup V_B|}p^{e(H)-e(H[V_R \cup V_G \cup V_O])}}{8n^{21}p^{45}} \left|  \sum_{\rho \in A_{R,G}} \sum_{\rho' \in A_{R,G,O}(\rho)} \left( \delta_{\rho'(\beta)\rho(\varsigma)}(\varphi_{t-3}) \pm n^{-\varepsilon/7}\delta_\infty(\varphi_{t-3}) \right) \right|\\
    & \leq \frac{n^{|V_P \cup V_B|}p^{e(H)-e(H[V_R \cup V_G \cup V_O])}}{8n^{21}p^{45}} \sum_{\rho \in A_{R,G}} \left| \sum_{\rho' \in A_{R,G,O}(\rho)} \delta_{\rho'(\beta)\rho(\varsigma)}(\varphi_{t-3}) \right|\\
    & \quad + \frac{n^{|V_P \cup V_B|}p^{e(H)-e(H[V_R \cup V_G \cup V_O])}}{8n^{21}p^{45}} |A_{R,G}| \max_{\rho \in A_{R,G}}|A_{R,G,O}(\rho)| n^{-\varepsilon/7} \delta_\infty(\varphi_{t-3})\\
    & \stackrel{\text{\cref{itm:color typicality}}}{\leq} \frac{n^{|V_P \cup V_B|}p^{e(H)-e(H[V_R \cup V_G \cup V_O])}}{8n^{21}p^{45}} \sum_{\rho \in A_{R,G}} \left| \sum_{\rho' \in A_{R,G,O}(\rho)} \delta_{\rho'(\beta)\rho(\varsigma)}(\varphi_{t-3}) \right|\\
    & \quad + np^2 n^{-\varepsilon/7} \delta_\infty(\varphi_{t-3}).
\end{align*}

We consider the first term. Let $\rho \in A_{R,G}$. Observe that each $\rho' \in A_{R,G,O}(\rho)$ is determined by the image of $\beta$ (i.e., by the image of the orange vertex). Every vertex in $N_G(\rho(\varsigma)) \setminus \rho(V_R \cup V_G)$ can play this role. Thus
\[
\sum_{\rho' \in A_{R,G,O}(\rho)} \delta_{\rho'(\beta)\rho(\varsigma)}(\varphi_{t-3}) = \sum_{x \in N_G(\rho(\varsigma))} \delta_{\rho(\varsigma)x}(\varphi_{t-3}) \pm |V_R \cup V_G| \delta_\infty(\varphi_{t-3}).
\]
Since $\varphi_{t-3}$ is vertex balanced we have $\sum_{x \in N_G(\rho(\varsigma))} \delta_{\rho(\varsigma)x}(\varphi_{t-3}) = 0$. Therefore
\[
\left| \sum_{\rho' \in A_{R,G,O}(\rho)} \delta_{\rho'(\beta)\rho(\varsigma)}(\varphi_{t-3}) \right| \leq |V_R \cup V_G| \delta_\infty(\varphi_{t-3}) \leq 7  \delta_\infty(\varphi_{t-3}).
\]
As a consequence
\begin{align*}
& \frac{n^{|V_P \cup V_B|}p^{e(H)-e(H[V_R \cup V_G \cup V_O])}}{8n^{21}p^{45}} \sum_{\rho \in A_{R,G}} \left| \sum_{\rho' \in A_{R,G,O}(\rho)} \delta_{\rho'(\beta)\rho(\varsigma)}(\varphi_{t-3}) \right|\\
& \leq \frac{n^{|V_P \cup V_B|}p^{e(H)-e(H[V_R \cup V_G \cup V_O])}}{8n^{21}p^{45}} |A_{R,G}| \times 7  \delta_\infty(\varphi_{t-3}) \stackrel{\text{\cref{itm:color typicality}}}{\leq} 7 p \delta_\infty(\varphi_{t-3}) \leq np^2 n^{-\varepsilon/7} \delta_\infty(\varphi_{t-3}).
\end{align*}

It follows that for $(i,j,m) \in \{(1,2,1),(1,7,7)\}$ we have
\[
|\delta(u,v;i,j,m)| \leq \gamma(u,v;i,j,m) + 9np^2n^{-\eps/7} \delta_\infty(\varphi_{t-3}) \leq 11np^2n^{-\eps/7} \delta_\infty(\varphi_{t-3}).
\]
By the induction hypothesis we have $\delta_\infty(\varphi_{t-3}) \leq 2^{-t+3}n^{-\varepsilon/8}$. Therefore for $(i,j,m) \in \{(1,2,1),\allowbreak(1,7,7)\}$ we have
\[
|\delta(u,v;i,j,m)| \leq \frac{1}{100} 2^{-t}np^2n^{-\eps/4}.
\]

It remains to handle $(i,j,m) \in \tripleLayerSet$ (see \cref{def:triple layer base set def}). The proof follows a similar path as that for $\doubleLayerConcentratedSet$.

Recall the definition of $\gamma$ from \eqref{eq:gamma def}. We will break down $\gamma(u,v;i,j,m)$ according to the image of $\vartheta_m \vartheta_{m-1}$. That is,
\begin{align*}
\gamma(u,v;i,j,m) & = \frac{1}{8n^{21}p^{45}} \left| \sum_{\alpha\beta \in E} \sum_{\ell=1}^{k} \sum_{\substack{\psi \in \chi_{i,j,\ell}:\\ \psi(\vartheta_m)\psi(\vartheta_{m-1})=\alpha\beta}} \delta_{\alpha\beta}(\varphi_{t-3}) \right|\\
& = \frac{1}{8n^{21}p^{45}} \left| \sum_{\alpha\beta \in E} \delta_{\alpha\beta}(\varphi_{t-3}) \sum_{\ell=1}^{k} |\{ \psi \in \chi_{i,j,\ell} : \psi(\vartheta_m)\psi(\vartheta_{m-1})=\alpha\beta \}|  \right|.
\end{align*}
We now observe that as long as $\alpha,\beta \notin \{u,v\}$ we have
\begin{align*}
\sum_{\ell=1}^{k} & |\{ \psi \in \chi_{i,j,\ell} : \psi(\vartheta_m)\psi(\vartheta_{m-1})=\alpha\beta \}|\\
& = X(W(i,j),G,\{d, a_7,\vartheta_m,\vartheta_{m-1}\},\{ d \mapsto v, a_7 \mapsto u, \vartheta_m \mapsto \alpha, \vartheta_{m-1} \mapsto \beta \})\\
& \quad + X(W(i,j),G,\{d, a_7,\vartheta_m,\vartheta_{m-1}\},\{ d \mapsto v, a_7 \mapsto u, \vartheta_m \mapsto \beta, \vartheta_{m-1} \mapsto \alpha \})\\
& \stackrel{\text{\cref{itm:12 triple-layer}}}{=} (1 \pm n^{-\varepsilon/6}) 2n^{20}p^{46}.
\end{align*}

Therefore
\begin{align*}
\gamma(u,v;i,j,m) & \leq \frac{1}{8n^{21}p^{45}} \left( \left| \sum_{e \in E} \delta_e(\varphi_{t-3}) (1 \pm n^{-\varepsilon/6}) \right| + (\deg_G(u) + \deg_G(v)) \delta_\infty(\varphi_{t-3}) \right) 2n^{20}p^{46}\\
& \stackrel{\text{\cref{itm:deg concentration}}}{\leq} \frac{p}{4n}\left| \sum_{e\in E}\delta_e(\varphi_{t-3}) \right| + n^{-\varepsilon/6} \frac{p}{4n} |E| \delta_\infty(\varphi_{t-3}) + \frac{p}{4n} 3np \delta_\infty(\varphi_{t-3}).
\end{align*}
Since the sum of discrepancies over all edges is zero, the first summand is zero. Hence
\[
\gamma(u,v;i,j,m) \leq \frac{p}{4n} \delta_\infty(\varphi_{t-3}) \left( n^{-\varepsilon/6} |E| + 3np \right) \stackrel{\text{\cref{itm:deg concentration}}}{\leq} \frac{p}{n} \delta_\infty(\varphi_{t-3}) n^{-\varepsilon/6} n^2p \leq np^2 n^{-\varepsilon/6} \delta_\infty(\varphi_{t-3}).
\]

It follows that for every $(i,j,m) \in \tripleLayerSet$ we have
\[
|\delta(u,v;i,j,m)| \leq \gamma(u,v;i,j,m) + 9np^2n^{-\eps/6} \delta_\infty(\varphi_{t-3}) \leq 10np^2n^{-\eps/6} \delta_\infty(\varphi_{t-3}).
\]
By the induction hypothesis we have $\delta_\infty(\varphi_{t-3}) \leq 2^{-t+3}n^{-\varepsilon/8}$. Therefore for $(i,j,m) \in \tripleLayerSet$
\[
|\delta(u,v;i,j,m)| \leq   \frac{1}{100} 2^{-t}np^2n^{-\eps/4}.
\]

As a consequence, by \eqref{eq:three layer overlap}, for $(i,j) \in \{(1,2),(1,7)\}$ we have
\[
|\delta(u,v;i,j)| \leq 2 \sum_{m=1}^7 |\delta(u,v;i,j,m)| + |R_{1,j}| \leq \frac{14}{100} 2^{-t}np^2n^{-\eps/4} + |R_{1,j}|.
\]

The upper bounds on $|\delta(u,v,i,j)|$ that we have proved for $(i,j) \in \doubleLayerConcentratedSet \cup \vertexBalanceBaseSet$ imply that for $(i,j) \in \{(1,1),(1,3),(1,4),(1,5),(1,6)\}$ we have
\[
|\delta(u,v;i,j)| \leq \frac{1}{100} 2^{-t}np^2n^{-\eps/4}.
\]
Therefore, by \eqref{eq:two layer overlap},
\begin{align*}
|\delta(u,v;1)| & \leq 2 \sum_{j=1}^7 |\delta(u,v;1,j)| + |R_1|\\
& \leq 2 \times \frac{1+14+1+1+1+1+14}{100} 2^{-t}np^2n^{-\eps/4} + |R_1| + 2\sum_{j\in\{2,7\}}|R_{1,j}|\\
& = \frac{66}{100} 2^{-t}np^2n^{-\eps/4} + |R_1| + 2\sum_{j\in\{2,7\}}|R_{1,j}|.
\end{align*}
Set $R_7=0$. In light of \eqref{eq:delta u v 7 bound} and \eqref{eq:delta most i}, we also have, for $i\in\{2,3,4,5,6,7\}$,
\[
|\delta(u,v;i)| \leq \frac{14}{100} 2^{-t}np^2n^{-\eps/4} + |R_i|.
\]
Therefore, by \eqref{eq:two layer remainder} and \eqref{eq:three layer remainder},
\begin{align*}
\left|\sum_{j=1}^k \delta_{uz_j}(\varphi_t) \right| & = 2 \left| \sum_{i=1}^7 \delta(u,v;i)(-1)^{i+1} \right| \leq 2 \sum_{i=1}^7 |\delta(u,v;i)|\\
& \leq 2 \times \frac{66+14\times6}{100} 2^{-t}np^2n^{-\eps/4} + 2\sum_{i=1}^6 |R_i| + 4\sum_{j\in\{2,7\}}|R_{1,j}|\\
& \leq \left(3+O(n^{-\eps/56})\right)2^{-t}np^2n^{-\eps/4} \leq 2^{-t}np^2n^{-\eps/7}.
\end{align*}
This completes our proof.
\end{proof}

\section{Concluding remarks}\label{sec:conclusion}

\subsection{Evidence in support of \cref{con:main conjecture}}

\begin{figure}
    \centering
    \includegraphics[width=\linewidth]{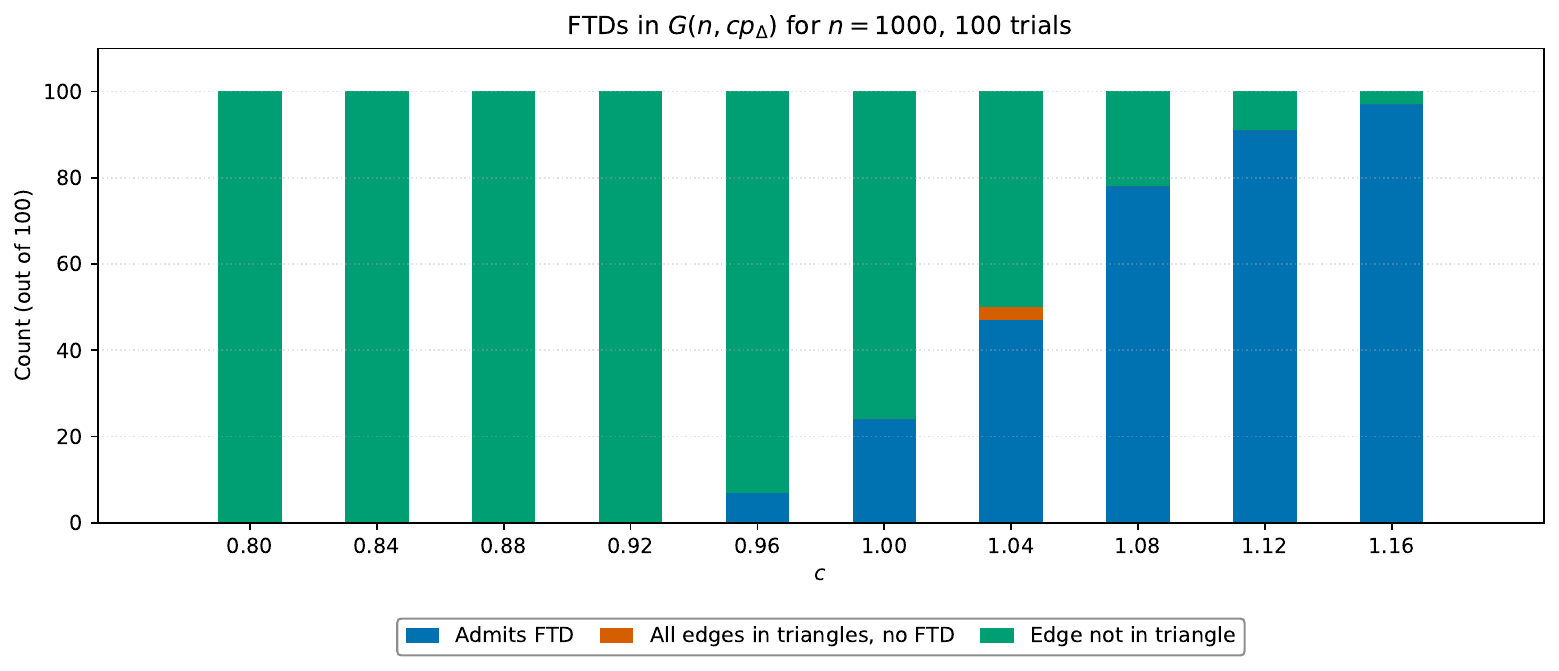}
    \caption{For each $c \in \{0.8,0.84,\ldots,1.16\}$, $100$ samples of $G(n, c\sqrt{3\log(n)/(2n)})$ were taken with $n=1000$. Each sample was tested for whether it supports an FTD and whether every edge is in a triangle. For all but $3$ samples, either there was an edge not in a triangle or the graphs admitted an FTD.}
    \label{fig:gnp experiments}
\end{figure}

Determining whether a graph admits an FTD can be formulated as a linear program, and as such can be solved efficiently both in theory and in practice. With this in mind, we tested \cref{con:main conjecture} as follows. We fixed $n=1000$. Let $p_\Delta \coloneqq \sqrt{3\log(n) / (2n)}$ be the conjectured threshold for $G(n,p)$ to admit an FTD. For each $c \in \{ 0.8, 0.84, \ldots, 0.96,1,1.04,\ldots,1.16 \}$ we sampled $100$ copies of $G(n,cp_\Delta)$ and tested whether they admit an FTD. For those that did not admit an FTD we tested whether this failure can be explained by the existence of an edge not contained in any triangle. The results of these experiments are presented in \cref{fig:gnp experiments}. Notably, among the $1000$ graphs considered, there were only three in which all edges were contained in triangles but did not admit FTDs (all three were found at $p=1.04 p_\Delta$). We view this as strong evidence in support of \cref{con:main conjecture}: empirically, w.h.p.\ random graphs either admit an FTD or have an edge not contained in a triangle.

We also tested a hitting time version of the conjecture: Given $n$, one can consider the random graph process in which edges are added in a uniformly random order until the first moment at which (the graph is nonempty and) every edge is contained in a triangle. For $n=1000$ we sampled $100$ graphs in this way. Of these, $99$ admitted FTDs. In light of this, we conjecture the following, which is stronger than \cref{con:main conjecture}.

\begin{conjecture}\label{con:hitting time}
    Let $n \in \N$ and let $G_0 \subseteq G_1 \subseteq \ldots \subseteq G_{\binom{n}{2}} = K_n$ be a random graph process in which $G_0$ is an empty graph on $n$ vertices and each $G_i$ is obtained from $G_{i-1}$ by adding a uniformly random edge. Let $\tau$ be the smallest $i>0$ in which every edge of $G_i$ is contained in a triangle. W.h.p., for every $t \geq \tau$, $G_t$ admits an FTD.
\end{conjecture}

\subsubsection*{Acknowledgement} The code for the experiments above was written with the assistance of Microsoft Copilot GPT-5. The linear programs were solved using Gurobi Optimizer \cite{gurobi}. The code can be found as part of the arXiv version of this paper \cite{mahabaduge2025fractional}.

\subsection{Future directions}

In our proof of \cref{thm:main} we described an algorithm to construct (a sequence of weightings converging to) an FTD. The algorithm used gadgets to shift weight between triangles. The gadgets we used (octagonal pinwheels) were specifically designed for fractional triangle decompositions in $G(n,n^{-4/11+o(1)})$. In this section we propose two directions to improve or generalize our approach.

The first suggestion, which we view as a plausible algorithm to find an FTD in $G(n,n^{-1/2+o(1)})$ (which would be optimal), is to simply replace the octagonal pinwheels with larger wheels. Specifically, fix $\varepsilon>0$ and let $G \sim G(n,n^{-1/2+\varepsilon})$. The first two steps of our algorithm (giving every triangle a uniform weight and then using bowties to correct the vertex defects) w.h.p.\ succeed in $G$. Unfortunately, for small $\varepsilon$, $G$ will be so sparse that many edges will not be contained in a copy of $W_8$. However, if $k \gg 1/\varepsilon$ is sufficiently large, then w.h.p.\ every edge is contained in many copies of $W_{2k}$. Thus, in the third phase of the algorithm, we can replace $W_8$ with $W_{2k}$. If under this operation the weightings converge to an FTD, it would follow that $G(n,n^{-1/2+o(1)})$ admits an FTD w.h.p. Unfortunately, analyzing the convergence of this algorithm seems daunting.

Our second suggestion might give a unified approach to finding fractional $K_k$-decompositions in random graphs, for any $k$. We first note that for a graph to admit a fractional (or integral) $K_k$-decomposition every edge must be contained in a copy of $K_k$. The threshold for this property in $G(n,p)$ is $n^{- (k-2)/(\binom{k}{2} - 1) + o(1)}$. This is a lower bound on the threshold for $G(n,p)$ to admit a fractional $K_k$-decomposition and it is natural to conjecture that this is indeed the threshold. 

We wonder whether the following simple algorithm can be used to find a sequence of weightings that converge to a fractional $K_k$-decomposition in $G \sim G(n,p)$, with $p$ above the conjectured threshold. Let $K$ be the number of (unlabeled) copies of $K_k$ in $G$. Let $\varphi_0 \equiv e(G) / (K \binom{k}{2})$. Note that since w.h.p. every edge of $G$ is in approximately $K \binom{k}{2} / e(G)$  copies of $K_k$, this is an approximate fractional $K_k$-decomposition. We now iteratively perform the following operation: For each edge $e$, subtract its discrepancy from all copies of $K_k$ that contain $e$, divided by the number of $K_k$s containing $e$. (In each iteration, this is done simultaneously for all edges.) In this way a sequence of weightings is produced. We wonder whether this sequence converges to a fractional $K_k$-decomposition of $G$.

Finally, we note that beyond determining the thresholds for $G(n,p)$ to admit a fractional or integral triangle decomposition, it is not even known that thresholds for these properties exist at all\footnote{This was true when the first version of this paper was posted to arXiv. See the next section for subsequent developments.}. Proving these facts, even without finding the locations of the thresholds, might be of substantial interest.

\subsubsection{Recent advances}

Subsequent to posting the first version of this paper to arXiv, there have been some exciting developments, which we relate here.

First, Joos and Smith \cite{joos2026fractional} proved that for every $\varepsilon>0$, w.h.p.\ $G(n,n^{-1/2+\varepsilon})$ admits an FTD. They furthermore determined the correct exponent for similar fractional clique decomposition problems in random hypergraphs. Their proof uses a ``clique reweighting'' algorithm that resembles the one proposed above.

Second, Tran \cite{tran2026threshold} proved that there is a constant $C$ such that for all $p \geq C\sqrt{\log n / n}$, w.h.p.\ $G(n,p)$ admits an FTD. Tran also proved analogous results for fractional clique decompositions of random hypergraphs. These determine the threshold for these properties up to multiplicative constants. In particular, thresholds for the properties exist. Tran's proof uses spectral techniques and is quite different from the approach taken in this paper.

\section*{AI acknowledgement}

The authors developed the proof and wrote the manuscript using traditional mathematical methods, without AI assistance. As described in the text, Microsoft Copilot (GPT-5) assisted with writing the code for the experiments. OpenAI Codex (GPT-6) assisted with implementing suggestions from the referee during revision, as well as checking mathematical arguments during the revision.

\bibliographystyle{amsplain}
\bibliography{biblio}

\end{document}